\documentclass[12pt]{article}
\usepackage{amsmath,amssymb,amsbsy}
\usepackage{graphicx}
\usepackage{multirow}
\usepackage{authblk}
\usepackage{algorithm}
\usepackage{xcolor}
\usepackage[normalem]{ulem}
\usepackage{tikz}
\usepackage{tkz-graph}
\usetikzlibrary{shapes.geometric}
\usepackage[letterpaper, left=2.5cm, right=2.5cm, top=2.5cm,bottom=2.5cm,dvips]{geometry}

\newtheorem{proposition}{Proposition}[section]
\newtheorem{theorem}[proposition]{Theorem}
\newtheorem{lemma}[proposition]{Lemma}

\newtheorem{conjecture}[proposition]{Conjecture}
\newtheorem{corollary}[proposition]{Corollary}

\newtheorem{problem}[proposition]{Problem}

\def\projn{$Pr_{n}(0,j)$}

\newcommand{\qed}{\hfill \rule{.1in}{.1in}}

\makeatletter
\def\imod#1{\allowbreak\mkern10mu({\operator@font mod}\,\,#1)}
\makeatother

\begin{document}

\title{On acyclic b-chromatic number of cubic graphs}

\author[1]{Marcin Anholcer\thanks{Partially supported by the National Science Centre of Poland under grant no. 2020/37/B/ST1/03298.}}
\author[2]{Sylwia Cichacz\thanks{The work of the author was supported by  the AGH University of Krakow under grant no. 16.16.420.054, funded by the Polish Ministry of Science and Higher Education.}}
\author[3]{Iztok Peterin\thanks{The author is partially supported by the Slovenian Research and Innovation Agency by program No. P1-0297 and is also with the Institute of Mathematics, Physics and Mechanics, Jadranska 19, 1000 Ljubljana, Slovenia.}}

\affil[1]{\scriptsize{}Pozna\'n University of Economics and Business, Institute of Informatics and Quantitative Economy}
\affil[ ]{Al.Niepodleg{\l}o\'sci 10, 61-875 Pozna\'n, Poland, \textit{marcin.anholcer@ue.poznan.pl}}
\affil[ ]{}
\affil[2]{AGH University, Faculty of Applied Mathematics}
\affil[ ]{Al. Mickiewicza 30, 30-059 Krak\'ow, Poland, \textit{cichacz@agh.edu.pl}}
\affil[ ]{}
\affil[3]{University of Maribor, Faculty of Electrical Engineering and Computer Science}
\affil[ ]{Koro\v{s}ka 46, 2000 Maribor, Slovenia,\textit{iztok.peterin@um.si}}

\maketitle

\begin{abstract}
Let $G$ be a graph. An acyclic $k$-coloring of $G$ is a map $c:V(G)\rightarrow \{1,\dots,k\}$ such that $c(u)\neq c(v)$ for any $uv\in E(G)$ and the subgraph induced by the vertices of any two colors $i,j\in \{1,\dots,k\}$ is a forest. If every vertex $v$ of a color class $V_i$ misses a color $\ell_v\in\{1,\dots,k\}$ in its closed neighborhood, then every $v\in V_i$ can be recolored with $\ell_v$ and we obtain a $(k-1)$-coloring of $G$. If a new coloring $c'$ is also acyclic, then such a recoloring is an acyclic recoloring step and $c'$ is in relation $\triangleleft_a$ with $c$. The acyclic b-chromatic number $A_b(G)$ of $G$ is the maximum number of colors in an acyclic coloring where no acyclic recoloring step is possible. Equivalently, it is the maximum number of colors in a minimum element of the transitive closure of $\triangleleft_a$. In this paper, we consider $A_b(G)$ of cubic graphs.
\end{abstract}

  \medskip


\section{Introduction}\label{sec_Intro}

A possible heuristic approach to properly color a graph is the one that starts with a trivial proper coloring, i.e., a coloring in which every vertex has its own color from $\{1,\dots,|V(G)|\}$. Then, given a color class $V_i$ of color $i$, a vertex $v\in V_i$ is called a b-vertex if it has all the colors other than $i$ in its neighborhood. If there is no b-vertex in $V_i$, then for every $v\in V_i$ there is at least one color $\ell_v$ other than $i$ not used in the neighbors of $v$, so each $v\in V_i$ can be recolored with $\ell_v$ to obtain a new proper coloring of $G$ using all colors except $i$. This removal of a color class is referred to as the recoloring step. Recoloring steps are performed until each color class has at least one b-vertex. The number of colors in the resulting coloring is an upper bound for $\chi(G)$. Its highest possible value is called the b-chromatic number $\varphi (G)$ of $G$. It was introduced by Irving and Manlove in 1999 \cite{irma-99}, who have shown, among others, that determining the b-chromatic number of a graph is an NP-hard problem. As demonstrated by Kratochv{\'{\i}}l, Tuza, and Voigt \cite{krtuvo-02}, the problem is still NP-hard for connected bipartite graphs. The most interesting results of the research on b-chromatic numbers for the present work are $d$-regular graphs. The main question for them is whether equality $\varphi(G)=d+1$ holds or not. It was proved in \cite{krtuvo-02} that this holds for all graphs on at least $d^4$ vertices. This was later improved to $2d^3$ vertices by Cabello and Jakovac \cite{caja-10} by the use of Hall's Marriage Theorem and later to $2d^3+2d-2d^2$ by El Sahili et al.\ \cite{ESKoMa}. In particular, this means that there exists only a finite number of $d$-regular graphs with $\varphi(G)<d+1$ and for $d=3$ there exists only four such graphs as shown by Jakovac and Klav\v zar \cite{JaKl}. This was later tested for small graphs and $d\in\{3,4,5,6,7\}$ by a hybrid evolutionary algorithm by Fister et al. \cite{FiPeMe}. The most famous conjecture on the topic was posted by El Sahili and Kouider \cite{ESKe} and later reformulated by Blidia et al. \cite{BlMaZe} (excluding the Petersen graph) and states that for every $d$-regular graph with girth at least five, except the Petersen graph, we have $\varphi(G)=d+1$. The latest improvement on it is by Detlaff et al. \cite{DFPR}.  We refer the reader to the survey \cite{JaPe} for additional topics related to b-coloring.

In the current paper, we consider an analogous problem, where the base is the acyclic coloring. In this case, one assumes not only the absence of monochromatic edges but also that every two color classes induce a forest. The minimum number of colors in such a coloring of graph $G$ is its \emph{acyclic chromatic number}, denoted by $A(G)$. Obviously, $A(G)\geq \chi(G)$ since every acyclic coloring is also a coloring of $G$. This kind of coloring was introduced by Gr\"unbaum \cite{ref_Gru}, who proved that $A(G)\leq 9$ for any planar graph $G$ and conjectured that $A(G)\leq 5$ in this case, which was finally proved by Borodin \cite{ref_Bor}. A conjecture attributed to Erd\"os (see \cite[p.89]{ref_JenTof}), stated that $A(G)=o(\Delta^2)$, where $\Delta=\Delta(G)$. It was proved by Alon, McDiarmid and Reed \cite{ref_AloMcDRee}, who showed that $A(G)\leq \lceil 50\Delta^{4/3}\rceil$. This result was recently improved to $A(G)\leq \frac{3}{2}\Delta^{4/3}+O(\Delta)$ by Gon\c{c}alves, Montassier and Pinlou \cite{ref_GonMonPin}. On the other hand, in \cite{ref_AloMcDRee} the existence of graphs satisfying $A(G)=\Omega(\Delta^{4/3}/(\log\Delta)^{1/3})$ has been proved. Other important and interesting results on acyclic colorings can be found in particular in the papers \cite{ref_MonNisWhiRah} ($A(G)$ of subdivisions of graphs and NP-completeness results), \cite{ref_AloMohSan}, where, in particular, the authors disproved a conjecture of Borodin (see \cite[p.70]{ref_JenTof}) about the equality $A(G)=\chi(G)$ for graphs embeddable on all surfaces other than a plane and \cite{ref_FerRas, ref_HocMon, ref_YadVarKotVen} (graphs with bounded degrees).

Motivated in particular by considerations about b-colorings, the authors introduced in \cite{AnCiPe} acyclic b-colorings, i.e., the colorings produced in the same way as b-colorings, with the extra condition that the recoloring step (called acyclic recoloring step) must produce an acyclic coloring. In particular, it was demonstrated that the problem is well-defined, i.e., for each finite graph $G$, every acyclic coloring of $G$ can be obtained by starting with a trivial coloring and performing a finite number of acyclic recoloring steps. The value of the respective graph invariant, the acyclic b-chromatic number $A_b(G)$, has been presented for several special graph families. Beside that, it was shown that the difference between $A_b(G)$ and the acyclic chromatic number $A(G)$, acyclic $m$-degree $m_a(G)$ (see $\cite{AnCiPe}$ for detailed definition), maximum degree $\Delta(G)$ and $\varphi (G)$, respectively, can be arbitrary big. We continue to develop the topic by analyzing the cubic graphs in this work. We believe that it is an important step on the way to obtaining more general results for regular graphs.

The paper is organized as follows. In the next section, we present basic notations and concepts; among others, we recall the necessary terms and results from \cite{AnCiPe}. In Section \ref{sec_cub} we show that all cubic graphs, besides one exception, have acyclic b-chromatic number $4$ or $5$. Moreover, in contrast to the b-chromatic number, there are infinitely many cubic graphs with $A_b(G)=4$. Two sections follow where we deal with generalized Petersen graphs and $(0,j)$-prisms, respectively. We conclude the paper with some final remarks, including open problems.


\section{Preliminaries}\label{sec_Preliminaries}

We deal only with finite and simple (without multiple edges and loops) graphs $G$ in this work, where $V(G)$ denotes the vertex set and $E(G)$ is the edge set of $G$. Let $n_G=\vert V(G)\vert$ and $m_G=\vert E(G)\vert$. For $v\in V(G)$ the \emph{open neighborhood} $N_G(v)$ equals to $\{u\in V(G):uv\in E(G)\}$ and the \emph{closed neighborhood} is $N_G[v]=N_G(v)\cup\{v\}$. The \emph{degree} of $v\in V(G)$ is $d_G(v)= \vert N_G(v)\vert$. By $\Delta(G)$ we denote the \emph{maximum degree} of a vertex from $V(G)$. If $d_G(v)=r$ for every $v\in V(G)$, then $G$ is an $r$-\emph{regular} graph. In particular, we will be interested here in 3-regular or \emph{cubic} graphs. For any $S\subseteq V(G)$, by $G[S]$ we denote the subgraph of $G$ induced by $S$. Set $\{1,\dots,k\}$ is denoted shortly by $[k]$.  

A (proper) \emph{vertex coloring} is a map $c:V(G)\rightarrow[k]$ where $c(x)\neq c(y)$ for every edge $xy\in E(G)$. The members of $[k]$ are called \emph{colors}. We color only vertices here and therefore omit the term "vertex" and call $c$ a coloring or a $k$-coloring of $G$. A \emph{trivial coloring} of $G$ is a coloring of $G$ with $n_G$ colors so that every vertex has its own color. The \emph{chromatic number} $\chi(G)$ of $G$ is the minimum integer $k$ for which there exists a $k$-coloring. The set $V_i=\{u\in V(G):c(u)=i\}$, for every $i\in [k]$, is called a \emph{color class} of $c$. Clearly, $V_1,\dots,V_k$ form a partition of $V(G)$ into independent sets. We will use the following notation: $V_{i,j}=V_i\cup V_j$ and $V_{i,j,\ell}=V_i\cup V_j\cup V_{\ell}$ for any $i,j,\ell\in [k]$. A $k$-coloring $c$ is an \emph{acyclic coloring} of $G$ if for any $i,j\in [k]$ graph $G[V_{i,j}]$ contains no cycles, thus $G[V_{i,j}]$ is a forest. Notice that $G[V_{i,i}]$ is a graph without edges. The \emph{acyclic chromatic number} $A(G)$ is the minimum number of colors of an acyclic coloring of $G$. Every acyclic coloring is also a coloring of $G$ and $A(G)\geq \chi(G)$ follows. 



Let $c$ be a $k$-coloring of a graph $G$. A vertex $v$ of $G$ is a b-\emph{vertex} for $c$, if all the colors are present in $N_G[v]$. If a vertex $v$ with $c(v)=i$ is not a b-vertex, then a color, say $\ell_v$, is missing in $N_G[v]$. If we recolor $v$ with $\ell_v$, then a slightly different coloring is obtained. This is possible for every vertex from $V_i$ whenever there is no b-vertex in $V_i$. Hence, if $V_i$ has no b-vertex for $c$, then $c_i:V(G)\rightarrow [k]\setminus\{i\}$ defined by 
\begin{equation}\label{recolor}
c_i(v)=\left\{ 
\begin{array}{ccc}
c(v) & \text{if } & c(v)\neq i, \\
\ell_v & \text{if } & c(v)=i, 
\end{array}
\right. 
\end{equation}
is a $(k-1)$-coloring of $G$. We refer to the above procedure as a \emph{recoloring step}. By the \emph{recoloring algorithm} we mean an iterative performing of recoloring steps while it is possible, where we start with a trivial coloring of $G$. The maximum number of colors obtained by the recoloring algorithm is called the \emph{b-chromatic number} of $G$, denoted by $\varphi (G)$. Clearly, as already observed in \cite{irma-99}, $\varphi (G)$ is the maximum number of colors in a coloring of $G$ where every color class contains a b-vertex. In contrast to $\varphi (G)$, $\chi (G)$ is the minimum number of colors obtained by the recoloring algorithm, and $\chi(G)\leq\varphi (G)$ follows. 

The approach to the b-chromatic number was restricted from colorings to acyclic colorings in \cite{AnCiPe}. Let $c:V(G)\rightarrow [k]$ be an acyclic coloring of $G$. If there exists a color class $V_i$, $i\in [k]$, such that for every vertex $v\in V_i$ there exists a color $\ell_v$, such that $c_i$ from (\ref{recolor}) is an acyclic coloring, then $c_i$ is obtained from $c$ by an \emph{acyclic recoloring step}. Similarly to before, the \emph{acyclic recoloring algorithm} consists of iteratively performing acyclic recoloring steps until it is possible, starting with a trivial coloring. The maximum number of colors obtained by the recoloring algorithm is called the \emph{acyclic b-chromatic number} of $G$ and denoted by $A_b(G)$, which is a kind of a dual to $A(G)$, which is the minimum number of colors obtained by the recoloring algorithm. Clearly, $A(G)\leq A_b(G)$.

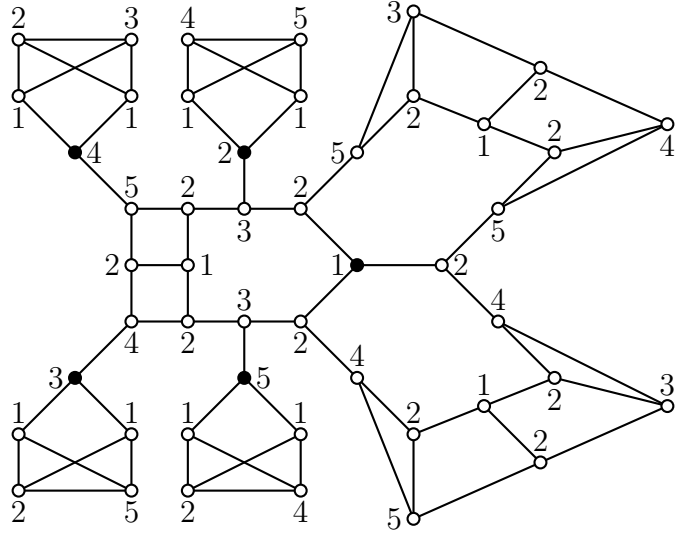
\begin{figure}[ht!]
\begin{center}
\begin{tikzpicture}[scale=0.75,style=thick,x=1cm,y=1cm]
\def\vr{3pt} 

\path (0,0) coordinate (a);
\path (1,0) coordinate (b);
\path (2,0) coordinate (c);
\path (3,0) coordinate (d);

\path (0,1) coordinate (a1);
\path (1,1) coordinate (b1);
\path (4,1) coordinate (c1);
\path (5.5,1) coordinate (d1);

\path (0,2) coordinate (a2);
\path (1,2) coordinate (b2);
\path (2,2) coordinate (c2);
\path (3,2) coordinate (d2);

\path (6.5,2) coordinate (a3);
\path (7.5,3) coordinate (b3);
\path (9.5,3.5) coordinate (c3);
\path (7.25,4.5) coordinate (d3);

\path (4,3) coordinate (a4);
\path (5,4) coordinate (b4);
\path (5,5.5) coordinate (c4);
\path (6.25,3.5) coordinate (d4);

\path (6.5,0) coordinate (a5);
\path (7.5,-1) coordinate (b5);
\path (9.5,-1.5) coordinate (c5);
\path (7.25,-2.5) coordinate (d5);

\path (4,-1) coordinate (a6);
\path (5,-2) coordinate (b6);
\path (5,-3.5) coordinate (c6);
\path (6.25,-1.5) coordinate (d6);

\path (-1,-1) coordinate (u);
\path (-2,-2) coordinate (v);
\path (0,-2) coordinate (w);
\path (-2,-3) coordinate (x);
\path (0,-3) coordinate (y);

\path (2,3) coordinate (u1);
\path (1,4) coordinate (v1);
\path (3,4) coordinate (w1);
\path (1,5) coordinate (x1);
\path (3,5) coordinate (y1);

\path (2,-1) coordinate (u2);
\path (1,-2) coordinate (v2);
\path (3,-2) coordinate (w2);
\path (1,-3) coordinate (x2);
\path (3,-3) coordinate (y2);

\path (-1,3) coordinate (u3);
\path (-2,4) coordinate (v3);
\path (0,4) coordinate (w3);
\path (-2,5) coordinate (x3);
\path (0,5) coordinate (y3);


\draw (b1) -- (a1) -- (a) -- (b) -- (c) -- (d) -- (c1) -- (d2) -- (c2) -- (b2) -- (a2) -- (a1);
\draw (b) -- (b1) -- (b2);
\draw (c1) -- (d1) -- (a3) -- (b3) -- (c3) -- (d3) -- (c4) -- (b4) -- (a4) -- (d2);
\draw (d3) -- (d4) -- (b4);
\draw (d4) -- (b3);
\draw (d1) -- (a5) -- (b5) -- (c5) -- (d5) -- (c6) -- (b6) -- (a6) -- (d);
\draw (d5) -- (d6) -- (b6);
\draw (d6) -- (b5);
\draw (c6) -- (a6);
\draw (c5) -- (a5);
\draw (a4) -- (c4);
\draw (a3) -- (c3);

\draw (a) -- (u) -- (v) -- (x) -- (y) -- (w) -- (u);
\draw (v) -- (y);
\draw (w) -- (x);
\draw (c2) -- (u1) -- (v1) -- (x1) -- (y1) -- (w1) -- (u1);
\draw (v1) -- (y1);
\draw (w1) -- (x1);
\draw (c) -- (u2) -- (v2) -- (x2) -- (y2) -- (w2) -- (u2);
\draw (v2) -- (y2);
\draw (w2) -- (x2);
\draw (a2) -- (u3) -- (v3) -- (x3) -- (y3) -- (w3) -- (u3);
\draw (v3) -- (y3);
\draw (w3) -- (x3);

\draw (a) [fill=white] circle (\vr);
\draw (b) [fill=white] circle (\vr);
\draw (c) [fill=white] circle (\vr);
\draw (d) [fill=white] circle (\vr);
\draw (a1) [fill=white] circle (\vr);
\draw (b1) [fill=white] circle (\vr);
\draw (c1) [fill=black] circle (\vr);
\draw (d1) [fill=white] circle (\vr);
\draw (a2) [fill=white] circle (\vr);
\draw (b2) [fill=white] circle (\vr);
\draw (c2) [fill=white] circle (\vr);
\draw (d2) [fill=white] circle (\vr);
\draw (a3) [fill=white] circle (\vr);
\draw (b3) [fill=white] circle (\vr);
\draw (c3) [fill=white] circle (\vr);
\draw (d3) [fill=white] circle (\vr);
\draw (a4) [fill=white] circle (\vr);
\draw (b4) [fill=white] circle (\vr);
\draw (c4) [fill=white] circle (\vr);
\draw (d4) [fill=white] circle (\vr);
\draw (a5) [fill=white] circle (\vr);
\draw (b5) [fill=white] circle (\vr);
\draw (c5) [fill=white] circle (\vr);
\draw (d5) [fill=white] circle (\vr);
\draw (a6) [fill=white] circle (\vr);
\draw (b6) [fill=white] circle (\vr);
\draw (c6) [fill=white] circle (\vr);
\draw (d6) [fill=white] circle (\vr);

\draw (u) [fill=black] circle (\vr);
\draw (v) [fill=white] circle (\vr);
\draw (x) [fill=white] circle (\vr);
\draw (y) [fill=white] circle (\vr);
\draw (w) [fill=white] circle (\vr);
\draw (u1) [fill=black] circle (\vr);
\draw (v1) [fill=white] circle (\vr);
\draw (x1) [fill=white] circle (\vr);
\draw (y1) [fill=white] circle (\vr);
\draw (w1) [fill=white] circle (\vr);
\draw (u2) [fill=black] circle (\vr);
\draw (x2) [fill=white] circle (\vr);
\draw (y2) [fill=white] circle (\vr);
\draw (w2) [fill=white] circle (\vr);
\draw (v2) [fill=white] circle (\vr);
\draw (u3) [fill=black] circle (\vr);
\draw (x3) [fill=white] circle (\vr);
\draw (y3) [fill=white] circle (\vr);
\draw (w3) [fill=white] circle (\vr);
\draw (v3) [fill=white] circle (\vr);

\draw[anchor = north] (a) node {$4$};
\draw[anchor = north] (b) node {$2$};
\draw[anchor = south] (c) node {$3$};
\draw[anchor = north] (d) node {$2$};
\draw[anchor = east] (a1) node {$2$};
\draw[anchor = west] (b1) node {$1$};
\draw[anchor = east] (c1) node {$1$};
\draw[anchor = west] (d1) node {$2$};
\draw[anchor = south] (a2) node {$5$};
\draw[anchor = south] (b2) node {$2$};
\draw[anchor = north] (c2) node {$3$};
\draw[anchor = south] (d2) node {$2$};
\draw[anchor = north] (a3) node {$5$};
\draw[anchor = south] (b3) node {$2$};
\draw[anchor = north] (c3) node {$4$};
\draw[anchor = north] (d3) node {$2$};
\draw[anchor = east] (a4) node {$5$};
\draw[anchor = north] (b4) node {$2$};
\draw[anchor = east] (c4) node {$3$};
\draw[anchor = north] (d4) node {$1$};
\draw[anchor = south] (a5) node {$4$};
\draw[anchor = north] (b5) node {$2$};
\draw[anchor = south] (c5) node {$3$};
\draw[anchor = south] (d5) node {$2$};
\draw[anchor = south] (a6) node {$4$};
\draw[anchor = south] (b6) node {$2$};
\draw[anchor = east] (c6) node {$5$};
\draw[anchor = south] (d6) node {$1$};

\draw[anchor = east] (u) node {$3$};
\draw[anchor = south] (v) node {$1$};
\draw[anchor = south] (w) node {$1$};
\draw[anchor = north] (x) node {$2$};
\draw[anchor = north] (y) node {$5$};
\draw[anchor = east] (u1) node {$2$};
\draw[anchor = north] (v1) node {$1$};
\draw[anchor = north] (w1) node {$1$};
\draw[anchor = south] (x1) node {$4$};
\draw[anchor = south] (y1) node {$5$};
\draw[anchor = west] (u2) node {$5$};
\draw[anchor = south] (v2) node {$1$};
\draw[anchor = south] (w2) node {$1$};
\draw[anchor = north] (x2) node {$2$};
\draw[anchor = north] (y2) node {$4$};
\draw[anchor = west] (u3) node {$4$};
\draw[anchor = north] (v3) node {$1$};
\draw[anchor = north] (w3) node {$1$};
\draw[anchor = south] (x3) node {$2$};
\draw[anchor = south] (y3) node {$3$};

\end{tikzpicture}
\end{center}
\caption{An acyclic b-coloring of a cubic graph. Black vertices are acyclic b-vertices: the one colored with $1$ is of type A and those colored with $2,3,4$ and $5$ are of type B.}
\label{less1}
\end{figure}

A vertex $v\in V(G)$ is an \emph{acyclic b-vertex} of color $c(v)=i$ if every possible recoloring of $v$ and possibly some other vertices of $V_i$ results in a bi-colored cycle which may contain $v$ or not (for the details see the definition of critical cycle systems in \cite{AnCiPe}). Clearly, we cannot recolor $v$ with any color that is present in $N_G(v)$. But it is possible to have a color outside of $N_G[v]$, and that recoloring step yields a bi-colored cycle. For this, observe the black vertices of Figure \ref{less1}. A black vertex of color $2$ has neighbors of color $1$ and $3$, and there are two cycles with colors $2-1-4-1$ and $2-1-5-1$, because of which $2$ cannot be recolored by $4$ and $5$, respectively, not to obtain a bi-colored cycle. Hence, the black vertex of color $2$ is an acyclic b-vertex for color $2$. It is similar to the black vertices of colors $3$, $4$, and $5$.

For the black vertex, call it $x$, of color $1$, there is only color $2$ in its neighborhood. For colors $3,4$ and $5$, there exist three cycles that contain $x$, each of which includes an extra vertex of color $1$. Moreover, this additional vertex can be recolored only by 3 ('left' cycle), only by $4$ ('bottom right' cycle), and only by $5$ ('upper right' cycle), respectively, if the coloring is supposed to remain acyclic (the other colors are blocked by the shorter cycles to which the vertices colored $1$ belong and by their neighborhoods). However, this would yield a bi-colored cycle if we recolor $x$ by $3$, $4$, or $5$. Hence, $x$ is an acyclic b-vertex for color $1$. Since there exist acyclic b-vertices of all colors for a coloring in Figure \ref{less1}, we cannot use an acyclic recoloring step anymore. For a more detailed description of acyclic b-vertices, we recommend Section 3 of \cite{AnCiPe}, where it was shown that $A_b(G)$ is the maximum integer $k$ such that there exists an acyclic $k$-coloring where every color class contains an acyclic b-vertex (Corollary 3.6).

While it is straightforward to see that $\varphi(G)\leq \Delta(G)+1$, $A_b(G)$ may be arbitrarily larger than $\Delta (G)$, see Theorem 4.2 from \cite{AnCiPe}. Nevertheless, the following general upper bound from \cite{AnCiPe} holds for every graph $G$:
$$A_b(G)\leq \frac{1}{2}(\Delta(G))^2+1.$$

In this work, we consider only cubic graphs. In such a case, the above upper bound transforms into 
\begin{equation}\label{bound}
A_b(G)\leq 5.
\end{equation}
Let $v$ be an acyclic b-vertex of a cubic graph $G$ that is not a b-vertex. Now, $v$ must have either one or two colors in $N_G(v)$ (otherwise no bi-colored cycle could appear). If all neighbors of $v$ have the same color, then we say that $v$ is an acyclic b-vertex of type A, and if there are two different colors in $N_G(v)$, then $v$ is an acyclic b-vertex of type B. See Figure \ref{less1} where the black vertex of color $1$ is of type A and the other black vertices are of type B. Further, let $v$ be an acyclic b-vertex of color $i$ such that color $j$ is not present in $N_G[v]$. A shortest cycle that blocks color $j$ is called a $j_v$-\emph{cycle}. A black vertex, say $v$, of color $5$ in Figure \ref{less1} has $2_v$- and $4_v$-cycles of length $4$.

 \begin{figure}[ht!]
\begin{center}
\begin{tikzpicture}[scale=0.7,style=thick,x=1cm,y=1cm]
\def\vr{3pt} 

\path (-5,0) coordinate (a); 
\path (-2,0) coordinate (b);
\path (-6,3) coordinate (c); 
\path (-1,3) coordinate (d);
\path (-3.5,5) coordinate (e);
\path (-4.25,1) coordinate (g); 
\path (-2.75,1) coordinate (h);
\path (-5,2.75) coordinate (i); 
\path (-2,2.75) coordinate (j);
\path (-3.5,3.75) coordinate (k);

\draw (a) -- (b); 
\draw (b) -- (d);
\draw (c) -- (e);
\draw (d) -- (e);
\draw (a) -- (c); 
\draw (a) -- (g);
\draw (b) -- (h);
\draw (c) -- (i);
\draw (d) -- (j);
\draw (e) -- (k);
\draw (k) -- (g);
\draw (j) -- (g);
\draw (h) -- (k);
\draw (i) -- (j);
\draw (i) -- (h);

\draw (a)  [fill=black] circle (\vr); \draw (b)  [fill=white] circle (\vr);
\draw (c)  [fill=black] circle (\vr); \draw (d)  [fill=white] circle (\vr);
\draw (e)  [fill=black] circle (\vr);
\draw (g)  [fill=white] circle (\vr); \draw (h)  [fill=white] circle (\vr);
\draw (i)  [fill=black] circle (\vr); \draw (j)  [fill=white] circle (\vr);
\draw (k)  [fill=white] circle (\vr);

\draw[anchor = east] (a) node {$3$};
\draw[anchor = west] (b) node {$4$};
\draw[anchor = north] (c) node {$2$};
\draw[anchor = north] (d) node {$2$};
\draw[anchor = west] (e) node {$1$};
\draw[anchor = west] (g) node {$2$};
\draw[anchor = east] (h) node {$2$};
\draw[anchor = south] (i) node {$4$};
\draw[anchor = south] (j) node {$3$};
\draw[anchor = west] (k) node {$3$};

\path (0,0) coordinate (a1);
\path (1.5,0) coordinate (b1);
\path (3,0) coordinate (c1);
\path (0,2) coordinate (d1);
\path (1.5,2) coordinate (e1);
\path (3,2) coordinate (f1);

\draw (a1) -- (b1) -- (c1) -- (f1) -- (e1) -- (d1) -- (a1);
\draw (b1) -- (e1);
\draw[-,line width=0.5pt] (a1) .. controls (1,-0.4) and (2,-0.4).. (c1);
\draw[-,line width=0.5pt] (d1) .. controls (1,2.4) and (2,2.4).. (f1);

\draw (a1) [fill=black] circle (\vr);
\draw (b1) [fill=black] circle (\vr);
\draw (c1) [fill=black] circle (\vr);
\draw (d1) [fill=white] circle (\vr);
\draw (e1) [fill=white] circle (\vr);
\draw (f1) [fill=white] circle (\vr);

\draw[anchor = north] (a1) node {$1$};
\draw[anchor = north] (b1) node {$2$};
\draw[anchor = north] (c1) node {$3$};
\draw[anchor = south] (d1) node {$3$};
\draw[anchor = south] (e1) node {$1$};
\draw[anchor = south] (f1) node {$2$};


\path (5,0) coordinate (a);
\path (6.5,0) coordinate (b);
\path (8,0) coordinate (c);
\path (5,2) coordinate (d);
\path (6.5,2) coordinate (e);
\path (8,2) coordinate (f);

\draw (a) -- (d) -- (c) -- (e) -- (a) -- (f) -- (b) -- (e);
\draw (b) -- (d);
\draw (c) -- (f);

\draw (a) [fill=black] circle (\vr);
\draw (b) [fill=black] circle (\vr);
\draw (c) [fill=black] circle (\vr);
\draw (d) [fill=black] circle (\vr);
\draw (e) [fill=white] circle (\vr);
\draw (f) [fill=white] circle (\vr);

\draw[anchor = north] (a) node {$1$};
\draw[anchor = north] (b) node {$2$};
\draw[anchor = north] (c) node {$3$};
\draw[anchor = south] (d) node {$4$};
\draw[anchor = south] (e) node {$4$};
\draw[anchor = south] (f) node {$4$};


\path (10,1) coordinate (g);
\path (10,3) coordinate (h);
\path (11.5,0) coordinate (i);
\path (11.5,4) coordinate (j);
\path (13,1) coordinate (k);
\path (13,3) coordinate (l);
\path (14.5,0) coordinate (m);
\path (14.5,4) coordinate (n);
\path (16,1) coordinate (o);
\path (16,3) coordinate (p);

\draw (g) -- (h) -- (j) -- (l) -- (n) -- (p) -- (j);
\draw (p) -- (o) -- (m) -- (k) -- (i) -- (g) -- (m);
\draw (l) -- (k);
\draw (g) -- (m);
\draw (i) -- (o);
\draw (h) -- (n);

\draw (g) [fill=black] circle (\vr);
\draw (h) [fill=white] circle (\vr);
\draw (i) [fill=black] circle (\vr);
\draw (j) [fill=white] circle (\vr);
\draw (k) [fill=black] circle (\vr);
\draw (l) [fill=white] circle (\vr);
\draw (m) [fill=white] circle (\vr);
\draw (n) [fill=white] circle (\vr);
\draw (o) [fill=black] circle (\vr);
\draw (p) [fill=white] circle (\vr);

\draw[anchor = east] (g) node {$1$};
\draw[anchor = east] (h) node {$4$};
\draw[anchor = north] (i) node {$2$};
\draw[anchor = south] (j) node {$2$};
\draw[anchor = west] (k) node {$3$};
\draw[anchor = west] (l) node {$1$};
\draw[anchor = north] (m) node {$2$};
\draw[anchor = south] (n) node {$2$};
\draw[anchor = west] (o) node {$4$};
\draw[anchor = west] (p) node {$3$};
\draw (12,1) node {$G_1$};

\end{tikzpicture}
\end{center}
\caption{Petersen graph, prism $K_2\Box K_3$, $K_{3,3}$ and $G_1$, and their acyclic b-colorings (acyclic b-vertices are black).}
\label{sporadic}
\end{figure}
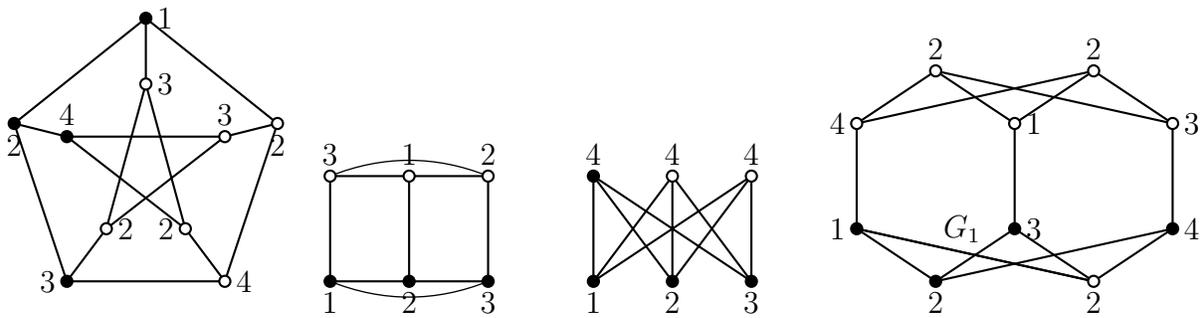


\section{All cubic graphs with one exception are $4$-acyclic b-colorable}\label{sec_cub}

In this section, we strongly rely on the following theorem and its proof by Jakovac and Klav\v zar \cite{JaKl}.

\begin{theorem}[\cite{JaKl}]\label{b-cubic}
For every cubic graph $G$ we have $\varphi(G)=4$, except for Petersen graph $P$, prism $K_2\Box K_3$, $K_{3,3}$ and $G_1$ from Figure \ref{sporadic}. Moreover, $\varphi(P)=\varphi(K_2\Box K_3)=\varphi(G_1)=3$ and $\varphi(K_{3,3})=2$. 
\end{theorem}

We start with a result about acyclic colorings of cubic graphs.  

\begin{theorem}\label{acycliccubic}
For every cubic graph $G$, we have $A(G)\leq 4$. 
\end{theorem}

\noindent {\textbf{Proof.}} Every cubic graph $G$ can be colored properly by four colors. Let $c:V(G)\rightarrow [4]$ be such a coloring. If $c$ is acyclic, then we are done. Otherwise, there exists a bi-colored cycle $C=v_1v_2\dots v_kv_1$, where $k$ is even. By a possible change of notation, we may assume that the vertices of $C_k$ are colored alternately with colors $1$ and $2$, where $c(v_1)=1$. Let $u_i$, $i\in [k]$, be the remaining neighbor of $v_i$. Suppose first that there exists $u_j$ with $c(u_j)\in\{1,2\}$, say $c(u_j)=1$ and let $\{a,b\}=\{3,4\}$. If $c(u_{j-1})=c(u_{j+1})=a$, then we define new coloring 
\begin{equation*}
c'(v)=\left\{ 
\begin{array}{ccc}
c(v) & \text{if } & v\neq v_j, \\
b & \text{if } & v=v_j, 
\end{array}
\right. 
\end{equation*}
and $C$ is not bi-colored by $c'$. Also any cycle that contains $u_{j-1}v_{j-1}v_jv_{j+1}u_{j+1}$ or $u_{j-1}v_{j-1}v_ju_j$ or $u_jv_jv_{j+1}u_{j+1}$ contains at least three colors $1,a,b$ and is not bi-colored by $c'$. So assume that $c(u_{j-1})=a$ and $c(u_{j+1})=b$. In such a case, we define coloring
\begin{equation*}
c''(v)=\left\{ 
\begin{array}{ccc}
c(v) & \text{if } & v\neq v_j, \\
b & \text{if } & v=v_{j-1}, 
\end{array}
\right. 
\end{equation*}
where $C$ is not bi-colored by $c''$. Moreover, no other new bi-colored cycles appear in $G$ by $c''$. 

Let now $c(u_i)\notin \{1,2\}$ for every $i\in [k]$ and suppose without loss of generality that $c(u_2)=3$. We define a new coloring 
\begin{equation*}
c'(v)=\left\{ 
\begin{array}{ccc}
c(v) & \text{if } & v\neq v_2, \\
4 & \text{if } & v=v_2, 
\end{array}
\right. 
\end{equation*}
and $C$ is not bi-colored by $c'$. On the other hand, $u_1v_1v_2v_3u_3$ can be a part of a new bi-colored cycle $C'$ when $c(u_1)=4=c(u_3)$. (Notice that this is not possible if $u_1=u_2$ or $u_3=U_2$.) In such a case, we define coloring
\begin{equation*}
c''(v)=\left\{ 
\begin{array}{ccc}
c'(v) & \text{if } & v\notin\{v_1,v_2,v_3\}, \\
3 & \text{if } & v\in\{v_1,v_3\}, \\
1 & \text{if } & v=v_2,
\end{array}
\right. 
\end{equation*}
where $C'$ is not bi-colored by $c''$. Again, no other bi-colored cycles are produced with either $c'$ or $c''$. 

Since the recoloring process defined above reduces the number of bi-colored cycles by one, the result now follows by induction on the number of bi-colored cycles in the initial $4$-coloring of $G$. \qed\\

With this result regarding acyclic colorings of cubic graphs and with Theorem \ref{b-cubic}, we can prove that $4\leq A_b(G)\leq 5$ for every cubic graph but the prism $K_2\Box K_3$.

\begin{theorem}\label{bacycliccubic}
For every cubic graph $G$ but the prism $K_2\Box K_3$, we have $A_b(G)\geq 4$. Moreover, $A_b(K_2\Box K_3)=3$. 
\end{theorem}

\noindent {\textbf{Proof.}} 
See Figure \ref{sporadic} for an acyclic b-coloring with four colors of Petersen graph, $K_{3,3}$ and $G_1$. For the prism $K_2\Box K_3$, assume that there exists an acyclic b-coloring of it with $4$ colors. We can also assume that one copy of $K_3$ is colored with colors $1$, $2$, and $3$ and that the neighbor of the vertex with color $1$ is colored by $4$. The two remaining vertices cannot be colored with $2$ and $3$, because this would result in a $2,3$-colored cycle. This means that color $1$ must be used on at least one of these vertices. Without loss of generality, let it be the remaining neighbor of the vertex colored $2$. The last remaining vertex must receive color $2$, because all other colors are already used in the neighborhood. But this way, there is no acyclic b-vertex, neither of color $3$ nor $4$, a contradiction. Hence, there is no acyclic b-coloring with four colors. It is trivial to see that there is no acyclic b-coloring of a prism with five colors. Therefore, we have $A_b(K_2\Box K_3)=3$ by the coloring presented in Figure \ref{sporadic}. 

\begin{figure}[ht!]
\begin{center}
\begin{tikzpicture}[scale=0.7,style=thick,x=1cm,y=1cm]
\def\vr{3pt} 

\path (0,5) coordinate (a1);
\path (1,5) coordinate (b1);
\path (2,5) coordinate (c1);
\path (3,5) coordinate (d1);
\path (1,4) coordinate (e1);
\path (2,4) coordinate (f1);
\path (3,4) coordinate (g1);
\path (1,3) coordinate (h1);
\path (2,3) coordinate (i1);
\path (3,3) coordinate (j1);
\path (2,2) coordinate (k1);
\path (3,2) coordinate (l1);

\draw (a1) -- (b1) -- (c1) -- (d1) -- (f1) -- (e1) -- (a1);
\draw (b1) -- (f1);
\draw (c1) -- (g1) --  (j1) -- (l1);
\draw (k1) -- (i1) -- (h1) -- (e1);
\draw (i1) -- (g1);
\draw[-,line width=0.5pt] (a1) .. controls (1,5.4) and (2,5.4).. (d1);
\draw[-,line width=0.5pt] (h1) .. controls (1.5,3.4) and (2.5,3.4).. (j1);

\draw (a1) [fill=white] circle (\vr);
\draw (b1) [fill=white] circle (\vr);
\draw (c1) [fill=white] circle (\vr);
\draw (d1) [fill=white] circle (\vr);
\draw (e1) [fill=white] circle (\vr);
\draw (f1) [fill=white] circle (\vr);
\draw (g1) [fill=black] circle (\vr);
\draw (h1) [fill=black] circle (\vr);
\draw (i1) [fill=black] circle (\vr);
\draw (j1) [fill=black] circle (\vr);
\draw (k1) [fill=white] circle (\vr);
\draw (l1) [fill=white] circle (\vr);

\draw[anchor = south] (a1) node {$2(1)$};
\draw (1,5.6) node {$4$};
\draw (2,5.6) node {$2$};
\draw[anchor = south] (d1) node {$3$};
\draw[anchor = east] (e1) node {$4$};
\draw[anchor = west] (f1) node {$1$};
\draw[anchor = west] (g1) node {$4$};
\draw[anchor = east] (h1) node {$2$};
\draw[anchor = west] (i1) node {$1$};
\draw[anchor = west] (j1) node {$3$};
\draw[anchor = east] (k1) node {$3$};
\draw[anchor = west] (l1) node {$1$};


\path (4,5) coordinate (a);
\path (5,5) coordinate (b);
\path (6,5) coordinate (c);
\path (7,5) coordinate (d);
\path (5,4) coordinate (e);
\path (6,4) coordinate (f);
\path (7,4) coordinate (g);
\path (7,3) coordinate (h);
\path (6,2) coordinate (i);
\path (7,2) coordinate (j);
\path (6,1) coordinate (k);
\path (7,1) coordinate (l);
\path (7,0) coordinate (m);
\path (7,-1) coordinate (n);
\path (6,-2) coordinate (o);
\path (7,-2) coordinate (p);

\draw (a) -- (b) -- (c) -- (d) -- (f) -- (e) -- (a);
\draw (b) -- (f);
\draw (c) -- (g) -- (h) -- (j) -- (l) -- (m) -- (n) -- (p);
\draw (n) -- (o);
\draw (m) -- (k) -- (i) -- (h);
\draw (i) -- (l);
\draw (j) -- (k);
\draw[-,line width=0.5pt] (a) .. controls (5,5.4) and (6,5.4).. (d);
\draw[-,line width=0.5pt] (e) .. controls (5.5,4.4) and (6.5,4.4).. (g);

\draw (a) [fill=white] circle (\vr);
\draw (b) [fill=white] circle (\vr);
\draw (c) [fill=white] circle (\vr);
\draw (d) [fill=white] circle (\vr);
\draw (e) [fill=white] circle (\vr);
\draw (f) [fill=white] circle (\vr);
\draw (g) [fill=black] circle (\vr);
\draw (h) [fill=black] circle (\vr);
\draw (i) [fill=white] circle (\vr);
\draw (j) [fill=white] circle (\vr);
\draw (k) [fill=white] circle (\vr);
\draw (l) [fill=white] circle (\vr);
\draw (m) [fill=black] circle (\vr);
\draw (n) [fill=black] circle (\vr);
\draw (o) [fill=white] circle (\vr);
\draw (p) [fill=white] circle (\vr);

\draw[anchor = south] (a) node {$4(3)$};
\draw(5,5.6) node {$1$};
\draw(6,5.6) node {$4$};
\draw[anchor = south] (d) node {$2$};
\draw[anchor = east] (e) node {$1$};
\draw[anchor = north] (f) node {$3$};
\draw[anchor = west] (g) node {$3$};
\draw[anchor = west] (h) node {$2$};
\draw[anchor = east] (i) node {$4$};
\draw[anchor = west] (j) node {$1$};
\draw[anchor = east] (k) node {$2$};
\draw[anchor = west] (l) node {$3$};
\draw[anchor = west] (m) node {$4$};
\draw[anchor = west] (n) node {$1$};
\draw[anchor = east] (o) node {$2$};
\draw[anchor = west] (p) node {$3$};


\path (8,5) coordinate (a2);
\path (9,5) coordinate (b2);
\path (10,5) coordinate (c2);
\path (11,5) coordinate (d2);
\path (9,4) coordinate (e2);
\path (10,4) coordinate (f2);
\path (11,4) coordinate (g2);
\path (11,3) coordinate (h2);
\path (10,2) coordinate (i2);
\path (11,2) coordinate (j2);
\path (9,1) coordinate (k2);
\path (10,1) coordinate (l2);
\path (11,1) coordinate (m2);

\draw (a2) -- (b2) -- (c2) -- (d2) -- (f2) -- (e2) -- (a2);
\draw (b2) -- (f2);
\draw (c2) -- (g2) -- (h2) -- (j2) -- (l2);
\draw (k2) -- (i2) -- (h2);
\draw (i2) -- (l2);
\draw (j2) -- (m2);
\draw[-,line width=0.5pt] (a2) .. controls (9,5.4) and (10,5.4).. (d2);
\draw[-,line width=0.5pt] (e2) .. controls (9.5,4.4) and (10.5,4.4).. (g2);

\draw (a2) [fill=white] circle (\vr);
\draw (b2) [fill=white] circle (\vr);
\draw (c2) [fill=white] circle (\vr);
\draw (d2) [fill=white] circle (\vr);
\draw (e2) [fill=white] circle (\vr);
\draw (f2) [fill=white] circle (\vr);
\draw (g2) [fill=black] circle (\vr);
\draw (h2) [fill=black] circle (\vr);
\draw (i2) [fill=black] circle (\vr);
\draw (j2) [fill=black] circle (\vr);
\draw (k2) [fill=white] circle (\vr);
\draw (l2) [fill=white] circle (\vr);
\draw (m2) [fill=white] circle (\vr);

\draw[anchor = south] (a2) node {$3(4)$};
\draw(9,5.6) node {$2$};
\draw(10,5.6) node {$3$};
\draw[anchor = south] (d2) node {$1$};
\draw[anchor = east] (e2) node {$1$};
\draw[anchor = north] (f2) node {$4$};
\draw[anchor = west] (g2) node {$2$};
\draw[anchor = west] (h2) node {$4$};
\draw[anchor = east] (i2) node {$3$};
\draw[anchor = west] (j2) node {$1$};
\draw[anchor = west] (k2) node {$1$};
\draw[anchor = west] (l2) node {$2$};
\draw[anchor = west] (m2) node {$3$};


\path (12,5) coordinate (a3);
\path (13,5) coordinate (b3);
\path (14,5) coordinate (c3);
\path (15,5) coordinate (d3);
\path (13,4) coordinate (e3);
\path (14,4) coordinate (f3);
\path (15,4) coordinate (g3);
\path (15,3) coordinate (h3);
\path (14,2) coordinate (i3);
\path (15,2) coordinate (j3);
\path (12,1) coordinate (k3);
\path (13,1) coordinate (l3);
\path (14,1) coordinate (m3);
\path (15,1) coordinate (n3);

\draw (a3) -- (b3) -- (c3) -- (d3) -- (f3) -- (e3) -- (a3);
\draw (b3) -- (f3);
\draw (c3) -- (g3) -- (h3) -- (j3) -- (n3);
\draw (k3) -- (i3) -- (h3);
\draw (i3) -- (l3);
\draw (j3) -- (m3);
\draw[-,line width=0.5pt] (a3) .. controls (13,5.4) and (14,5.4).. (d3);
\draw[-,line width=0.5pt] (e3) .. controls (13.5,4.4) and (14.5,4.4).. (g3);

\draw (a3) [fill=white] circle (\vr);
\draw (b3) [fill=white] circle (\vr);
\draw (c3) [fill=white] circle (\vr);
\draw (d3) [fill=white] circle (\vr);
\draw (e3) [fill=black] circle (\vr);
\draw (f3) [fill=white] circle (\vr);
\draw (g3) [fill=black] circle (\vr);
\draw (h3) [fill=black] circle (\vr);
\draw (i3) [fill=black] circle (\vr);
\draw (j3) [fill=white] circle (\vr);
\draw (k3) [fill=white] circle (\vr);
\draw (l3) [fill=white] circle (\vr);
\draw (m3) [fill=white] circle (\vr);
\draw (n3) [fill=white] circle (\vr);

\draw[anchor = south] (a3) node {$4$};
\draw(13,5.6) node {$1$};
\draw(14,5.6) node {$3$};
\draw[anchor = south] (d3) node {$2(1)$};
\draw[anchor = east] (e3) node {$1$};
\draw[anchor = north] (f3) node {$3$};
\draw[anchor = west] (g3) node {$2$};
\draw[anchor = west] (h3) node {$4$};
\draw[anchor = east] (i3) node {$3$};
\draw[anchor = west] (j3) node {$1$};
\draw[anchor = west] (k3) node {$1$};
\draw[anchor = west] (l3) node {$2$};
\draw[anchor = west] (m3) node {$3$};
\draw[anchor = west] (n3) node {$4$};


\path (16.5,5) coordinate (a4);
\path (17.5,5) coordinate (b4);
\path (18.5,5) coordinate (c4);
\path (19.5,5) coordinate (d4);
\path (17.5,4) coordinate (e4);
\path (18.5,4) coordinate (f4);
\path (19.5,4) coordinate (g4);
\path (17.5,3) coordinate (h4);
\path (18.5,3) coordinate (i4);
\path (19.5,3) coordinate (j4);
\path (18.5,2) coordinate (k4);
\path (19.5,2) coordinate (l4);
\path (18.5,1) coordinate (m4);
\path (19.5,1) coordinate (n4);

\draw (a4) -- (b4) -- (c4) -- (d4) -- (f4); 
\draw (e4) -- (a4);
\draw (b4) -- (f4) -- (i4);
\draw (c4) -- (g4) --  (j4) -- (e4) -- (h4) -- (g4);
\draw (h4) -- (k4) -- (i4);
\draw (i4) -- (l4) -- (m4);
\draw (j4) -- (k4);
\draw (l4) -- (n4);
\draw[-,line width=0.5pt] (a4) .. controls (17,5.4) and (18,5.4).. (d4);

\draw (a4) [fill=white] circle (\vr);
\draw (b4) [fill=white] circle (\vr);
\draw (c4) [fill=white] circle (\vr);
\draw (d4) [fill=white] circle (\vr);
\draw (e4) [fill=white] circle (\vr);
\draw (f4) [fill=black] circle (\vr);
\draw (g4) [fill=white] circle (\vr);
\draw (h4) [fill=white] circle (\vr);
\draw (i4) [fill=black] circle (\vr);
\draw (j4) [fill=white] circle (\vr);
\draw (k4) [fill=black] circle (\vr);
\draw (l4) [fill=black] circle (\vr);
\draw (m4) [fill=white] circle (\vr);
\draw (n4) [fill=white] circle (\vr);

\draw[anchor = south] (a4) node {$3(2)$};
\draw(17.5,5.6) node {$1$};
\draw(18.5,5.6) node {$3$};
\draw[anchor = south] (d4) node {$2$};
\draw[anchor = east] (e4) node {$4$};
\draw[anchor = south] (f4) node {$3$};
\draw[anchor = west] (g4) node {$2$};
\draw[anchor =east] (h4) node {$3$};
\draw[anchor = east] (i4) node {$4$};
\draw[anchor = west] (j4) node {$1$};
\draw[anchor = west] (k4) node {$2$};
\draw[anchor = west] (l4) node {$1$};
\draw[anchor = west] (m4) node {$2$};
\draw[anchor = west] (n4) node {$3$};

\end{tikzpicture}
\end{center}
\caption{Corrected b-colorings of \cite{JaKl}: first graph of the second line of Figure 14 and first, fourth, and fifth graphs from Figure 15; that are now acyclic. On the last graph, we have corrected a minor coloring error in the third graph from Figure 13 from \cite{JaKl}. The original colors are indicated in brackets, and black vertices represent b-vertices.}
\label{recoloring}
\end{figure}
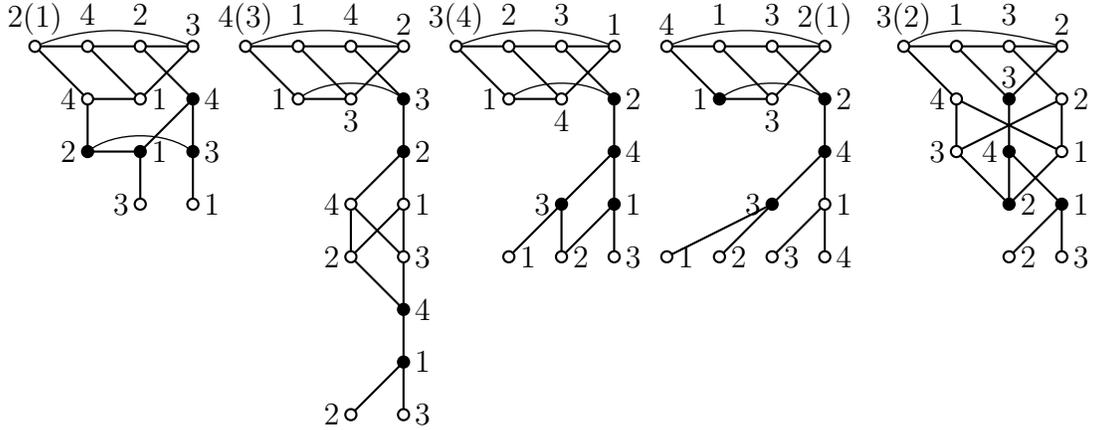

In the remainder, we assume that $G$ is none of the graphs mentioned before. We will strongly rely on the proof of Theorem \ref{b-cubic} presented in \cite{JaKl}. This proof is organized as follows. One starts to color a cubic graph on a shortest cycle $C$ and then continues to color vertices close to $C$ so that every color has a b-vertex. The remainder of the graph is then colored in an arbitrary manner, such as a greedy algorithm. In \cite{JaKl}, one can find a large number of pictures of graphs, and for each of them, besides the four exceptions, one can find a b-coloring in the first (small) part. So, we only need to be careful about bi-colored cycles that appear in figures from \cite{JaKl}. Luckily, there are only four such graphs in all figures, and they are the first four graphs in Figure \ref{recoloring}. The original coloring is in brackets, and it is easy to see the bi-colored cycles (they are always four-cycles). The new color (outside of the brackets) shows that it is easy to recolor these four graphs in an acyclic fashion. Notice that black vertices represent b-vertices of appropriate colors. Since the rest of the graph is colored using a greedy algorithm, one can also produce some bi-colored cycles. For such cycles (if they exist), we can use the recoloring from the proof of Theorem \ref{acycliccubic} and cancel all of them. Hence, we have constructed an acyclic b-coloring of $G$ with four colors, which implies $A_b(G)\geq 4$. 
\qed\\

Next, we show that the behavior of $A_b(G)$ differs from the behavior of $\varphi(G)$ on cubic graphs. In contrast to Theorem \ref{b-cubic}, we can find infinitely many cubic graphs with $A_b(G)=4$. For this, we need the following lemma.

\begin{figure}[ht!]
\begin{center}
\begin{tikzpicture}[scale=0.5,style=thick,x=1cm,y=1cm]
\def\vr{3pt} 

\path (-9,1) coordinate (u3);
\path (-9,-1) coordinate (v3);
\path (-7.5,1) coordinate (x3);
\path (-7.5,-1) coordinate (y3);
\path (-6,1) coordinate (w3);
\path (-6,-1) coordinate (z3);
\path (-10,0) coordinate (f3);

\path (-0.5,0) coordinate (a);
\path (-2,0) coordinate (b);
\path (-3,1) coordinate (c);
\path (-3,-1) coordinate (d);

\path (2.5,2) coordinate (a1);
\path (2.5,-2) coordinate (b1);
\path (3.5,1) coordinate (c1);
\path (3.5,-1) coordinate (d1);
\path (4.5,0) coordinate (e1);
\path (5.5,0) coordinate (f1);

\path (10.5,2) coordinate (a2);
\path (10.5,-2) coordinate (b2);
\path (11.5,1) coordinate (c2);
\path (11.5,-1) coordinate (d2);
\path (12.5,0) coordinate (e2);
\path (13.5,0) coordinate (f2);

\path (10,2.5) coordinate (u);
\path (10,1.5) coordinate (v);
\path (9,2.5) coordinate (x);
\path (9,1.5) coordinate (y);
\path (8,2.5) coordinate (w);
\path (8,1.5) coordinate (z);

\path (10,-2.5) coordinate (u1);
\path (10,-1.5) coordinate (v1);
\path (9,-2.5) coordinate (x1);
\path (9,-1.5) coordinate (y1);
\path (8,-2.5) coordinate (w1);
\path (8,-1.5) coordinate (z1);

\path (14,0.5) coordinate (u2);
\path (14,-0.5) coordinate (v2);
\path (15,0.5) coordinate (x2);
\path (15,-0.5) coordinate (y2);
\path (16,0.5) coordinate (w2);
\path (16,-0.5) coordinate (z2);
\draw [->] (0.5,0) -- (1.5,0);
\draw [->] (6,0) -- (7,0);

\draw (a) -- (b) -- (c);
\draw (b) -- (d);
\draw (a1) -- (c1) -- (e1) -- (f1);
\draw (b1) -- (d1) -- (e1);
\draw (c1) -- (d1);
\draw (a2) -- (c2) -- (e2) -- (f2);
\draw (b2) -- (d2) -- (e2);
\draw (c2) -- (d2);

\draw (a2) -- (u) -- (x) -- (w) -- (z) -- (y) -- (v) -- (a2);
\draw (u) -- (v);
\draw (x) -- (z);
\draw (y) -- (w);

\draw (b2) -- (u1) -- (x1) -- (w1) -- (z1) -- (y1) -- (v1) -- (b2);
\draw (u1) -- (v1);
\draw (x1) -- (z1);
\draw (y1) -- (w1);

\draw (f2) -- (u2) -- (x2) -- (w2) -- (z2) -- (y2) -- (v2) -- (f2);
\draw (u2) -- (v2);
\draw (x2) -- (z2);
\draw (y2) -- (w2);

\draw (f3) -- (u3) -- (x3) -- (w3) -- (z3) -- (y3) -- (v3) -- (f3);
\draw (u3) -- (v3);
\draw (x3) -- (z3);
\draw (y3) -- (w3);

\draw (a) [fill=white] circle (\vr);
\draw (b) [fill=white] circle (\vr);
\draw (c) [fill=white] circle (\vr);
\draw (d) [fill=white] circle (\vr);
\draw (a1) [fill=white] circle (\vr);
\draw (b1) [fill=white] circle (\vr);
\draw (c1) [fill=white] circle (\vr);
\draw (d1) [fill=white] circle (\vr);
\draw (e1) [fill=white] circle (\vr);
\draw (f1) [fill=white] circle (\vr);
\draw (a2) [fill=white] circle (\vr);
\draw (b2) [fill=white] circle (\vr);
\draw (c2) [fill=white] circle (\vr);
\draw (d2) [fill=white] circle (\vr);
\draw (e2) [fill=white] circle (\vr);
\draw (f2) [fill=white] circle (\vr);
\draw (u) [fill=white] circle (\vr);
\draw (v) [fill=white] circle (\vr);
\draw (x) [fill=white] circle (\vr);
\draw (y) [fill=white] circle (\vr);
\draw (z) [fill=white] circle (\vr);
\draw (w) [fill=white] circle (\vr);
\draw (u1) [fill=white] circle (\vr);
\draw (v1) [fill=white] circle (\vr);
\draw (x1) [fill=white] circle (\vr);
\draw (y1) [fill=white] circle (\vr);
\draw (z1) [fill=white] circle (\vr);
\draw (w1) [fill=white] circle (\vr);
\draw (u2) [fill=white] circle (\vr);
\draw (x2) [fill=white] circle (\vr);
\draw (y2) [fill=white] circle (\vr);
\draw (z2) [fill=white] circle (\vr);
\draw (w2) [fill=white] circle (\vr);
\draw (v2) [fill=white] circle (\vr);
\draw (u3) [fill=white] circle (\vr);
\draw (x3) [fill=white] circle (\vr);
\draw (y3) [fill=white] circle (\vr);
\draw (z3) [fill=white] circle (\vr);
\draw (w3) [fill=white] circle (\vr);
\draw (v3) [fill=white] circle (\vr);
\draw (f3) [fill=white] circle (\vr);

\draw (14,-1.5) node {$C(T)$};
\draw (-2,-1.5) node {$T$};
\draw[anchor = north] (a) node {$x$};
\draw[anchor = east] (b) node {$v$};
\draw[anchor = east] (c) node {$y$};
\draw[anchor = east] (d) node {$z$};
\draw[anchor = east] (a1) node {$y$};
\draw[anchor = east] (b1) node {$z$};
\draw[anchor = east] (c1) node {$b$};
\draw[anchor = east] (d1) node {$c$};
\draw[anchor = north] (e1) node {$a$};
\draw[anchor = north] (f1) node {$x$};
\draw[anchor = east] (f3) node {$w$};
\draw (-8,0) node {$H_3$};
\draw[anchor = south] (u3) node {$4$};
\draw[anchor = north] (v3) node {$3$};
\draw[anchor = south] (x3) node {$3$};
\draw[anchor = north] (y3) node {$4$};
\draw[anchor = south] (w3) node {$1$};
\draw[anchor = north] (z3) node {$2$};

\end{tikzpicture}
\end{center}
\caption{Graph $H_3$ and the construction of a graph $C(T)$ from a cubic tree $T\cong K_{1,3}$.}
\label{less}
\end{figure}
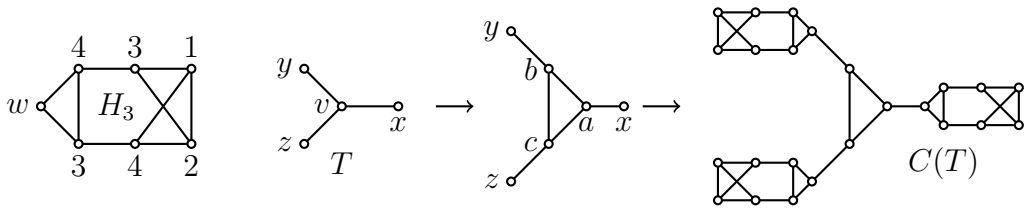

\begin{lemma}\label{simple}
Let $G$ be a cubic graph. If there exists a subgraph of $G$ that is isomorphic to $H_3$ from Figure \ref{less}, where $w$ is the vertex of degree two in $H_3$, then $w$ is not an acyclic b-vertex in an acyclic b-coloring with $5$ colors of $G$ (if it exists). 
\end{lemma}

\noindent {\textbf{Proof.}} Vertex $w$ can have at most four colors in its closed neighborhood. To be an acyclic b-vertex in an acyclic b-coloring with five colors, it must be an acyclic b-vertex of type A or type B. Clearly, $w$ is a cut-vertex of $G$. Hence, $w$ must be type B, where the neighbors of $w$ in $H_3$ have the same color, a contradiction because the mentioned neighbors are adjacent. So, $w$ is not an acyclic b-vertex in a coloring with five colors. 
\qed\\

A tree $T\ncong K_1$ is called \textit{cubic} if every inner vertex of $T$ is of degree three. This means that a cubic tree contains only inner vertices of degree three and leaves (of degree one). The smallest cubic tree is $K_2$, that is, without inner vertices, and next is $K_{1,3}$ with one inner vertex. 

From a cubic tree $T$, we construct the cubic graph $C(T)$ as follows. First, we replace every inner vertex $v$ having neighbors $x,y,z$ by a triangle $abc$, where edges $ax,by$ and $cz$ are added between the triangle and $N_T(v)$. Then, we join a copy of $H_3$ with every leaf of $T$ by identifying this leaf with the respective copy of $w$. As an example, we presented the construction of $C(K_{1,3})$ in Figure \ref{less}. 

\begin{theorem}\label{notupper}
If $T$ is a cubic tree, then $A_b(C(T))=4$. 
\end{theorem}

\noindent {\textbf{Proof.}} 
Let $T$ be a cubic tree. First, we present an acyclic b-coloring of $C(T)$ with four colors. We color every copy of $H_3$ as it is presented in Figure \ref{less}. The remaining vertices can be colored greedily because $C(T)$ is a cubic graph and we have four colors. This is a proper coloring with b-vertices of every color in every copy of $H_3$. Moreover, this coloring is also acyclic. This is easily seen in every copy of $H_3$, and all the remaining cycles in $C(T)$ are $3$-cycles. So, $A_b(C(T))\geq 4$. 

By (\ref{bound}) we know that $A_b(C(T))\leq 5$. Let us assume that $A_b(C(T))=5$. Every color class needs to have an acyclic b-vertex that is not a b-vertex. Clearly, this is not a vertex on a triangle that originates in an inner vertex of $T$, because every such vertex is contained in only one (odd) cycle in $C(T)$. Moreover, also vertex $w$ from any copy of $H_3$ in Figure \ref{less} is not an acyclic b-vertex in a 5-coloring by Lemma \ref{simple}. Any other vertex $u$ of any copy of $H_3$ cannot be an acyclic b-vertex, since no bi-colored cycle can contain $w$ and there is no other induced even-length cycle in $H_3$. Hence, $A_b(C(T))<5$, which leads to the desired equality. \qed\\

A direct consequence of Theorem \ref{notupper} is the following corollary, which reveals another contrast between the b-chromatic number and the acyclic b-chromatic number. Recall that there are precisely four cubic graphs for which $\chi_b(G)<4$ by Theorem \ref{b-cubic}. 

\begin{corollary}\label{infinite}
The number of cubic graphs $G$ with $A_b(G)<5$ is not finite. 
\end{corollary}



\section{Generalized Petersen Graphs}

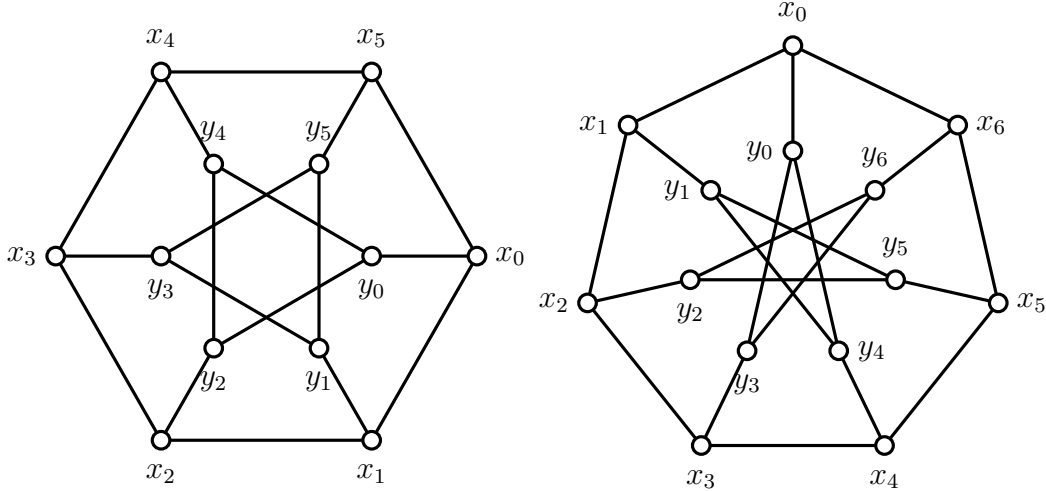
\begin{figure}[ht!]
\begin{center}
\begin{tikzpicture}
  [scale=.7,auto=left,every node/.style={shape=circle,minimum size = 1pt, very thick}]
  \node[draw=black,scale=.6] (x0) at (4,0) {};
  \node[draw=black,scale=.6] (x1) at (2,-3.5) {};
  \node[draw=black,scale=.6] (x2) at (-2,-3.5) {};
  \node[draw=black,scale=.6] (x3) at (-4,0) {};
  \node[draw=black,scale=.6] (x4) at (-2,3.5) {};
  \node[draw=black,scale=.6] (x5) at (2,3.5) {};
  \node[draw=black,scale=.6] (y0) at (2,0) {};
  \node[draw=black,scale=.6] (y1) at (1,-1.75) {};
  \node[draw=black,scale=.6] (y2) at (-1,-1.75) {};
  \node[draw=black,scale=.6] (y3) at (-2,0) {};
  \node[draw=black,scale=.6] (y4) at (-1,1.75) {};
  \node[draw=black,scale=.6] (y5) at (1,1.75) {};
  \draw[anchor = west] (x0) node {$x_0$};
  \draw[anchor = north] (x1) node {$x_1$};
  \draw[anchor = north] (x2) node {$x_2$};
  \draw[anchor = east] (x3) node {$x_3$};
  \draw[anchor = south] (x4) node {$x_4$};
  \draw[anchor = south] (x5) node {$x_5$};
  \draw[anchor = north] (y0) node {$y_0$};
  \draw[anchor = north] (y1) node {$y_1$};
  \draw[anchor = north] (y2) node {$y_2$};
  \draw[anchor = north] (y3) node {$y_3$};
  \draw[anchor = south] (y4) node {$y_4$};
  \draw[anchor = south] (y5) node {$y_5$};
  
  \foreach \from/\to in {x0/x1,x1/x2,x2/x3,x3/x4,x4/x5,x5/x0,
  y0/y2,y1/y3,y2/y4,y3/y5,y4/y0,y5/y1,
  x0/y0,x1/y1,x2/y2,x3/y3,x4/y4,x5/y5}
  \draw[-,very thick] (\from) -- (\to);
  \node[draw=black,scale=.6] (z0) at (10,4) {};
  \node[draw=black,scale=.6] (z1) at (6.87,2.49) {};
  \node[draw=black,scale=.6] (z2) at (6.1,-0.89) {};
  \node[draw=black,scale=.6] (z3) at (8.26,-3.6) {};
  \node[draw=black,scale=.6] (z4) at (11.74,-3.6) {};
  \node[draw=black,scale=.6] (z5) at (13.9,-0.89) {};
  \node[draw=black,scale=.6] (z6) at (13.13,2.49) {};
  \node[draw=black,scale=.6] (w0) at (10,2) {};
  \node[draw=black,scale=.6] (w1) at (8.44,1.25) {};
  \node[draw=black,scale=.6] (w2) at (8.05,-0.45) {};
  \node[draw=black,scale=.6] (w3) at (9.13,-1.8) {};
  \node[draw=black,scale=.6] (w4) at (10.87,-1.8) {};
  \node[draw=black,scale=.6] (w5) at (11.95,-0.45) {};
  \node[draw=black,scale=.6] (w6) at (11.56,1.25) {};
  \draw[anchor = south] (z0) node {$x_0$};
  \draw[anchor = east] (z1) node {$x_1$};
  \draw[anchor = east] (z2) node {$x_2$};
  \draw[anchor = north] (z3) node {$x_3$};
  \draw[anchor = north] (z4) node {$x_4$};
  \draw[anchor = west] (z5) node {$x_5$};
  \draw[anchor = west] (z6) node {$x_6$};
  \draw[anchor = east] (w0) node {$y_0$};
  \draw[anchor = east] (w1) node {$y_1$};
  \draw[anchor = north] (w2) node {$y_2$};
  \draw[anchor = north] (w3) node {$y_3$};
  \draw[anchor = west] (w4) node {$y_4$};
  \draw[anchor = south] (w5) node {$y_5$};
  \draw[anchor = south] (w6) node {$y_6$};
  
  \foreach \from/\to in {z0/z1,z1/z2,z2/z3,z3/z4,z4/z5,z5/z6,z6/z0,
  w0/w3,w1/w4,w2/w5,w3/w6,w4/w0,w5/w1,w6/w2,
  z0/w0,z1/w1,z2/w2,z3/w3,z4/w4,z5/w5,z6/w6}
    \draw[-,very thick] (\from) -- (\to);
  
\end{tikzpicture}
\end{center}
\caption{Graphs $G(6,2)$ and $G(7,3)$.}
\label{generalizedPetersen}
\end{figure}

The generalized Petersen graphs $G(n,k)$, where $1\leq k<n/2$, introduced by Coxeter in \cite{ref_Cox} and named by Watkins \cite{ref_Wat}, are the graphs on $2n$ vertices $\{x_0,\dots, x_{n-1},y_0,\dots,y_{n-1}\}$, where the edge set consists of the polygon $\{x_ix_{i+1}:0\leq i\leq n-1\}$, the star polygon $\{y_iy_{i+k}:0\leq i\leq n-1\}$ and the spokes $\{x_iy_i:0\leq i\leq n-1\}$, where the sums are taken modulo $n$. In this notation, the ordinary Petersen graph is denoted as $G(5,2)$, see the left graph of Figure \ref{sporadic}. The examples of $G(6,2)$ and $G(7,3)$ are given in Figure \ref{generalizedPetersen}.

From (\ref{bound}) it follows that $A_b(G(n,k))\leq 5$, however,
it is easy to prove that an acyclic b-coloring of $G(5,2)$ with $5$ colors does not exist. Indeed, if we have an acyclic b-vertex $u$ of type A in $G(5,2)$, then the color of $x$ is the only appearance of this color, and the three colors not in the neighborhood of $x$ have no acyclic b-vertex. So, assume that $x$ is a type B acyclic b-vertex with two neighbors $y_1$ and $y_2$ in one color and the third neighbor $z$ in another. There exist two vertex-disjoint cycles $xy_1abcy_2x$ and $xy_1defy_2x$ of length $6$ because the girth of $G(5,2)$ is $5$ and there is not enough space for longer cycles. This implies that $z$ is the only vertex in its color class. Therefore, it cannot be an acyclic b-vertex, a contradiction. An acyclic b-coloring of $G(5,2)$ with four colors is presented in Figure \ref{sporadic}, see the left graph. 

However, if $n$ is large enough with respect to $k\geq 3$, one can prove that there exists an acyclic b-coloring of $G(n,k)$ with $5$ colors.

\begin{theorem}\label{thm_petersen}
If $k\geq 3$ and $n\geq 5(2k+(-1)^k)$, then $A_b(G(n,k))=5$.
\end{theorem}

\noindent {\textbf{Proof.}} 
It is enough to construct the respective acyclic b-colorings with five colors. Note that throughout the proof, the addition and subtraction are taken modulo $n$ if necessary. No vertex can be a b-vertex, and we construct an acyclic b-coloring where every acyclic b-vertex $x_j$ is of type B.

\begin{figure}[ht!]
\begin{center}

\resizebox{\textwidth}{!}{
\begin{tikzpicture}
  [scale=0.7,auto=left,every node/.style={shape=circle,minimum size = 1pt, very thick}]
  \node[draw=black,scale=.6] (xj) at (0,0) {};
  \node[draw=black,scale=.6] (xj-1) at (-3,0) {};
  \node[draw=black,scale=.6] (xj-k+1) at (-7,0) {};
  \node[draw=black,scale=.6] (xj-k) at (-10,0) {};
  \node[draw=black,scale=.6,label=0:$y_{j-k}:c_j^1$] (yj-k) at (-10,3) {};
  \node[draw=black,scale=.6,label=0:$y_{j}:c_j^3$] (yj) at (0,3) {};

  \node[draw=black,scale=.6] (xj+1) at (3,0) {};

  \node[draw=black,scale=.6,label=180:$y_{j-1}:c_j^2$] (yj-1) at (-3,3) {};
  \node[draw=black,scale=.6] (xj+k-1) at (7,0) {};
  \node[draw=black,scale=.6] (xj+k) at (10,0) {};
  \node[draw=black,scale=.6,label=180:$y_{j+k-1}:c_j^3$] (yj+k-1) at (7,3) {};
  \node[draw=black,scale=.6,label=180:$y_{j+k}:c_j^2$] (yj+k) at (10,3) {};
  
  \draw (7,-.6) node {$x_{j+k-1}:c_j^2$};
  \draw (10,-.6) node {$x_{j+k}:c_j^3$};
  \draw (3,-.6) node {$x_{j+1}:c_j^4$};
  \draw (0,-.6) node {$x_j:c_j^0$};
  \draw (-3,-.6) node {$x_{j-1}:c_j^3$};
  \draw (-7,-.6) node {$x_{j-k+1}:c_j^1$};
  \draw (-10,-.6) node {$x_{j-k}:c_j^3$};
 
  \foreach \from/\to in {xj/xj-1,xj-k+1/xj-k,xj-k/yj-k,xj/yj,xj-1/yj-1,xj/xj+1,xj+k-1/xj+k,xj+k-1/yj+k-1,xj+k/yj+k}
  \draw[-,very thick] (\from) -- (\to);

   \foreach \from/\to in {xj-1/xj-k+1}
  \draw[dashed,very thick] (\from) -- (\to);

  \foreach \from/\to in {yj-k/yj,yj-1/yj+k-1,yj/yj+k}
  \draw[-,very thick] (\from) [out = 45, in = 135] to (\to);
  
\end{tikzpicture}
}

\end{center}
\caption{An acyclic b-vertex $x_j$ of $G(n,k)$ for odd $k$.}
\label{genPetersenOdd}
\end{figure}
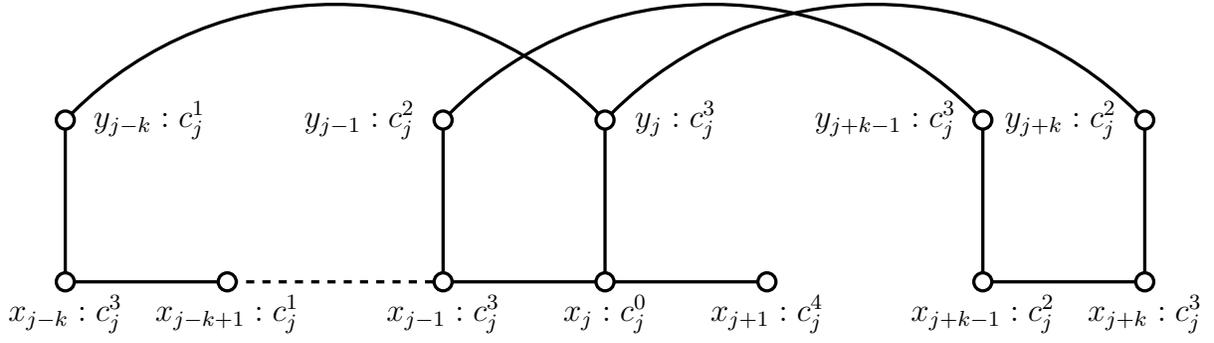

If $k$ is odd, then consider the following cycles $C_j^1=x_jx_{j-1}\dots x_{j-k+1}x_{j-k}y_{j-k}y_jx_j$ and $C_j^2=x_jx_{j-1}y_{j-1}y_{j+k-1}x_{j+k-1}x_{j+k}y_{j+k}y_jx_j$, where the colors assigned to the vertices are $c(x_j)=c_j^0$, $c(y_j)=c(y_{j+k-1})=c(x_{j+k})=c(x_{j-i})=c_j^3$ for odd $1\leq i\leq k$, $c(y_{j-k})=c(x_{j-i})=c_j^1$ for even $2\leq i\leq k-1$, $c(y_{j-1})=c(x_{j+k-1})=c(y_{j+k})=c_j^2$ and $c(x_{j+1})=c_j^4$, where $c_j^0$, $c_j^1$, $c_j^2$, $c_j^3$ and $c_j^4$ are pairwise distinct. Now, $x_j$ is an acyclic b-vertex for color $c_j^0$, since the colors $c_j^3$ and $c_j^4$ are in its neighborhood, while $C_j^1$ is $(c_j^1)_{x_j}$-cycle and $C_j^2$ is $(c_j^2)_{x_j}$-cycle. The situation is illustrated in Figure \ref{genPetersenOdd}.

One can repeat the presented configuration by choosing the acyclic b-vertices to be $x_{j_i}$, where $j_i=(2k-1)i$ for $0\leq i\leq 4$. This identifies the vertices $x_{j_i+k-1}$, $y_{j_i+k-1}$ and $x_{j_i+k}$ with $x_{j_{i+1}-k}$, $y_{j_{i+1}-k}$ and $x_{j_{i+1}-k+1}$, respectively. This can be done, since the set of identities $(c_{j_i}^2,c_{j_i}^3)= (c_{j_{i+1}}^3,c_{j_{i+1}}^1)$ defined for $0\leq i\leq 4$ can be easily extended to the equalities of permutations giving the correct colorings of consecutive segments of a graph even if one cyclically joins them (i.e., when $n=10k-5$), e.g. as: 
\begin{eqnarray*}
(c_{j_i}^0,c_{j_i}^1,c_{j_i}^2,c_{j_i}^3,c_{j_i}^4)= 
(c_{j_{i+1}}^4,c_{j_{i+1}}^0,c_{j_{i+1}}^3,c_{j_{i+1}}^1,c_{j_{i+1}}^2), 0\leq i\leq 4.
\end{eqnarray*}

In the case $n=10k-4$, i.e., when the first and the last segments meet on exactly one pair of vertices (and so the respective mapping in one case must satisfy $(c_{j_4}^3,c_{j_4}^2)= (c_{j_{0}}^3,c_{j_{0}}^1)$ instead of $(c_{j_4}^2,c_{j_4}^3)= (c_{j_{0}}^3,c_{j_{0}}^1)$), it is enough to redefine the above equalities for $i\in\{2,4\}$:
\begin{eqnarray*}
(c_{j_i}^0,c_{j_i}^1,c_{j_i}^2,c_{j_i}^3,c_{j_i}^4)=
\begin{cases}
			(c_{j_{i+1}}^4,c_{j_{i+1}}^0,c_{j_{i+1}}^3,c_{j_{i+1}}^1,c_{j_{i+1}}^2), & \text{if $i\in\{0,1,3\}$,}\\
            (c_{j_{i+1}}^2,c_{j_{i+1}}^0,c_{j_{i+1}}^3,c_{j_{i+1}}^1,c_{j_{i+1}}^4), & \text{if $i=2$,}\\
            (c_{j_{i+1}}^4,c_{j_{i+1}}^0,c_{j_{i+1}}^1,c_{j_{i+1}}^3,c_{j_{i+1}}^2), & \text{if $i=4$.}\\
		 \end{cases}
\end{eqnarray*}
Note that in both cases the colors $c_{j_i}^0$ are distinct in every permutation, so there are indeed acyclic b-vertices for all the colors in the presented colorings. It follows that finding an appropriate coloring is possible if $n\geq 10k-5$.

\begin{figure}[ht!]
\begin{center}

\resizebox{\textwidth}{!}{
\begin{tikzpicture}
  [scale=0.6,auto=left,every node/.style={shape=circle,minimum size = 1pt, very thick}]
  \node[draw=black,scale=.6] (xj) at (0,0) {};
  \node[draw=black,scale=.6] (xj-1) at (-3,0) {};
  \node[draw=black,scale=.6] (yj-1) at (-3,3) {};
  \node[draw=black,scale=.6] (yj-k-1) at (-13,3) {};
  \node[draw=black,scale=.6] (xj-k-1) at (-13,0) {};
  \node[draw=black,scale=.6] (xj-k) at (-10,0) {};
  \node[draw=black,scale=.6] (yj-k) at (-10,3) {};
  \node[draw=black,scale=.6] (yj) at (0,3) {};

  \node[draw=black,scale=.6] (xj+1) at (3,0) {};

  \node[draw=black,scale=.6] (xj-2) at (-6,0) {};
  \node[draw=black,scale=.6] (yj-2) at (-6,3) {};
  \node[draw=black,scale=.6] (yj+k-2) at (7,3) {};
  \node[draw=black,scale=.6] (xj+k-2) at (7,0) {};
  \node[draw=black,scale=.6] (xj+k-1) at (10,0) {};
  \node[draw=black,scale=.6] (xj+k) at (13,0) {};  
  \node[draw=black,scale=.6] (yj+k) at (13,3) {};
  \draw (-13.8,4) node {$y_{j-k-1}:c_j^3$};
  \draw (-2.2,3.6) node {$y_{j-1}:c_j^1$};
  \draw (-8.3,2.6) node {$y_{j-k}:c_j^1$};
  \draw (-7,3.7) node {$y_{j-2}:c_j^3$};
  \draw[anchor = west] (yj) node {$y_{j}:c_j^3$};
  \draw[anchor = east] (yj+k) node {$y_{j+k}:c_j^2$};
  \draw[anchor = east] (yj+k-2) node {$y_{j+k-2}:c_j^2$};
  \draw (0,-.6) node {$x_j:c_j^0$};
  \draw (-3,-.6) node {$x_{j-1}:c_j^3$};
  \draw (-13.5,-.6) node {$x_{j-k-1}:c_j^1$};
  \draw (-9.5,-.6) node {$x_{j-k}:c_j^3$};
  \draw (3,-.6) node {$x_{j+1}:c_j^4$};
  \draw (-6,-.6) node {$x_{j-2}:c_j^2$};
  \draw (6.5,-.6) node {$x_{j+k-2}:c_j^3$};
   \draw (10,-.6) node {$x_{j+k-1}:c_j^2$};
  \draw (13.3,-.6) node {$x_{j+k}:c_j^3$};
    
  \foreach \from/\to in {xj/xj-1,xj-1/yj-1,yj-k-1/xj-k-1,xj-k-1/xj-k,xj-k/yj-k,yj/xj,xj/xj+1,
  xj-1/xj-2,xj-2/yj-2,yj+k-2/xj+k-2,xj+k-2/xj+k-1,xj+k-1/xj+k,xj+k/yj+k}
  \draw[-,very thick] (\from) -- (\to);


  \foreach \from/\to in {yj-k/yj,yj-k-1/yj-1,yj-2/yj+k-2,yj/yj+k}
  \draw[-,very thick] (\from) [out = 45, in = 135] to (\to);
  
\end{tikzpicture}
}

\end{center}
\caption{An acyclic b-vertex $x_j$ of $G(n,k)$ for even $k$.}
\label{genPetersenEven}
\end{figure}
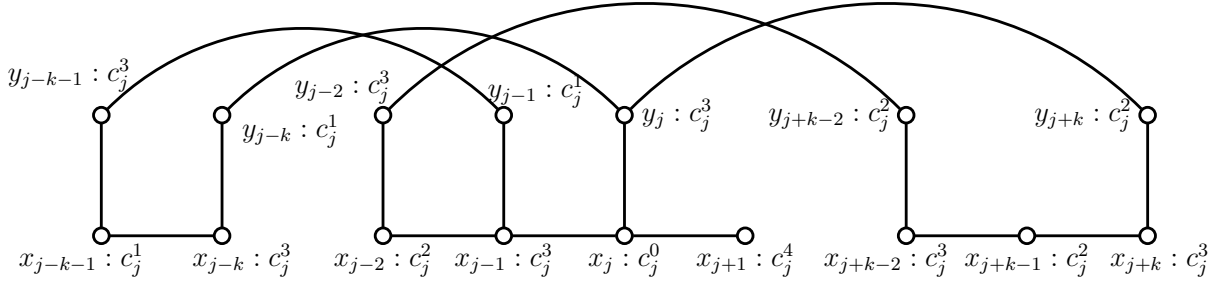

If $k$ is even, then we consider cycles $C_j^1=x_jx_{j-1}y_{j-1}y_{j-k-1}x_{j-k-1} x_{j-k}y_{j-k}y_jx_j$ and $C_j^2=x_jx_{j-1}x_{j-2}y_{j-2}y_{j+k-2}x_{j+k-2}x_{j+k-1}x_{j+k}y_{j+k}y_jx_j$, where the colors assigned to the vertices are $c(x_j)=c_j^0$, $c(x_{j-1})=c(y_{j-k-1})=c(x_{j-k})=c(y_j)=c(y_{j-2})=c(x_{j+k-2})=c(x_{j+k})=c_j^3$, $c(y_{j-1})=c(x_{j-k-1})=c(y_{j-k})=c_j^1$, $c(x_{j-2})=c(y_{j+k-2})=c(x_{j+k-1})=c(y_{j+k})=c_j^2$ and $c(x_{j+1})=c_j^4$, where $c_j^0$, $c_j^1$, $c_j^2$, $c_j^3$ and $c_j^4$ are pairwise distinct. Then $x_j$ is an acyclic b-vertex for color $c_j^0$, since it has neighbors colored $c_j^3$ and $c_j^4$, while $C_j^1$ is $(c_j^1)_{x_j}$-cycle and $C_j^2$ is $(c_j^2)_{x_j}$-cycle. The situation is illustrated in Figure \ref{genPetersenEven}.

This time, one can repeat the presented configuration by choosing the acyclic b-vertices to be $x_{j_i}$, where $j_i=(2k+1)i$ for $0\leq i\leq 4$. This identifies the vertices $x_{j_i+k}$ and $y_{j_i+k}$ with $x_{j_{i+1}-k-1}$ and $y_{j_{i+1}-k-1}$, respectively, which can be done, since the set of mappings $(c_{j_i}^2,c_{j_i}^3)\rightarrow (c_{j_{i+1}}^3,c_{j_{i+1}}^1)$ defined for $0\leq i\leq 4$ can be obviously extended to the permutations giving the correct colorings of consecutive segments of graph even if one cyclically joins them. Actually, one can use the same mappings as in the case of odd $k$. Note that also in this case, the colors $c_{j_i}^0$ are distinct in every permutation, so there exists an acyclic b-vertex for each color in the presented coloring. Thus, finding an appropriate coloring is possible if $n\geq 10k+5$.

\qed\\

A special kind of generalized Petersen graph is a prism, that is, $G(n,1)$. In this case, it can be proven that any acyclic b-coloring must use at most four colors. It follows that in fact $A_b(G(n,1))=4$ with only one exception. Indeed, since $G(3,1)=K_2\Box K_3$, $A_b(G(3,1))=3$ by Theorem \ref{bacycliccubic}. All the remaining cases are covered by the following result.

\begin{theorem}
For any integer $n\geq 4$ we have $A_b(G(n,1))=4$.\label{thm_prism}
\end{theorem}

\noindent {\textbf{Proof.}} 
By Theorem \ref{bacycliccubic}, it is enough to show that $A_b(G(n,1))\leq 4$. We already know that $A_b(G(n,1))\leq 5$ by (\ref{bound}), so it is enough to show that the value $5$ is impossible.

Suppose that $A_b(G(n,1))=5$ for some $n\geq 4$ and let $c$ be an appropriate acyclic b-coloring with five colors. Assume without loss of generality that for some $0\leq j\leq n-1$, $x_j$ is an acyclic b-vertex for color $c_1$. It cannot be a b-vertex since it has only $3$ neighbors. This means that there are only three different cases to consider (the remaining ones follow by symmetry).\\

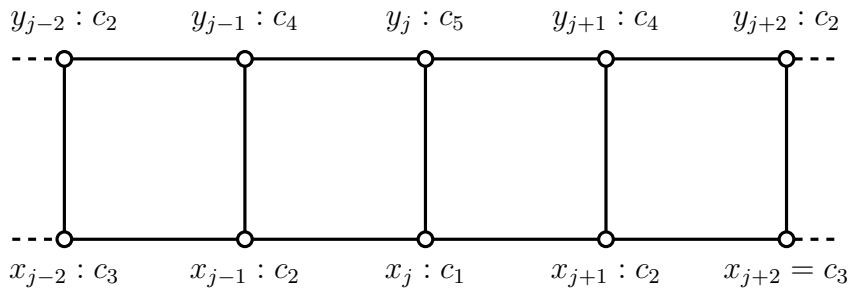
\begin{figure}[ht!]
\begin{center}

\begin{tikzpicture}
  [scale=0.8,auto=left,every node/.style={shape=circle,minimum size = 1pt, very thick}]
  \node[draw=black,scale=.5] (xj) at (0,0) {};
  \node[draw=black,scale=.5] (xj-1) at (-3,0) {};
  \node[draw=black,scale=.5] (xj-2) at (-6,0) {};
  \node[draw=white,scale=.5] (xj-3) at (-7,0) {};
  \node[draw=black,scale=.5] (xj+1) at (3,0) {};
  \node[draw=black,scale=.5] (xj+2) at (6,0) {};
  \node[draw=white,scale=.5] (xj+3) at (7,0) {};
  
  \node[draw=black,scale=.5] (yj) at (0,3) {};
  \node[draw=black,scale=.5] (yj-1) at (-3,3) {};
  \node[draw=black,scale=.5] (yj-2) at (-6,3) {};
  \node[draw=white,scale=.5] (yj-3) at (-7,3) {};
  \node[draw=black,scale=.5] (yj+1) at (3,3) {};
  \node[draw=black,scale=.5] (yj+2) at (6,3) {};
  \node[draw=white,scale=.5] (yj+3) at (7,3) {};

  \draw (0,-.6) node {$x_j:c_1$};
  \draw (-3,-.6) node {$x_{j-1}:c_2$};
  \draw (-6,-.6) node {$x_{j-2}:c_3$};
  \draw (3,-.6) node {$x_{j+1}:c_2$};
  \draw (6,-.6) node {$x_{j+2}=c_3$};

   \draw (0,3.6) node {$y_j:c_5$};
  \draw (-3,3.6) node {$y_{j-1}:c_4$};
  \draw (-6,3.6) node {$y_{j-2}:c_2$};
  \draw (3,3.6) node {$y_{j+1}:c_4$};
  \draw (6,3.6) node {$y_{j+2}:c_2$};    
  \foreach \from/\to in {xj-2/xj-1,xj-1/xj,xj/xj+1,xj+1/xj+2,yj-2/yj-1,yj-1/yj,yj/yj+1,yj+1/yj+2,xj-2/yj-2,xj-1/yj-1,xj/yj,xj+1/yj+1,xj+2/yj+2}
  \draw[-,very thick] (\from) -- (\to);

   \foreach \from/\to in {xj-3/xj-2,yj-3/yj-2,xj+2/xj+3,yj+2/yj+3}
  \draw[dashed,very thick] (\from) -- (\to);

  
\end{tikzpicture}

\end{center}
\caption{Prisms: $c(x_{j-1})=c_2=c(x_{j+1})\neq c(y_j)=c_5$.}
\label{prismCase2}
\end{figure}

\noindent{}\textbf{Case 1.} $x_j$ is type B vertex with $c(x_{j-1})=c(x_{j+1})\neq c(y_{j})$ (see Figure \ref{prismCase2}).\\
Let $c(x_{j-1})=c(x_{j+1})=c_2$ and $c(y_{j})=c_5$. There must exist two internally disjoint paths with common endpoints $x_{j-1}$ and $x_{j+1}$, one colored with $c_2$ and $c_3$ and the other one with $c_2$ and $c_4$. The only possible configuration is presented in Figure \ref{prismCase2}. Note that this time the internal vertices of the path colored with $c_2$ and $c_3$ are precisely all the vertices $x_i$, where $i\notin\{j-1,j,j+1\}$, while the internal vertices of the path colored with $c_2$ and $c_4$ are precisely all the vertices $y_i$, where $i\neq j$ (and $n$ must be even in this case). So, all the vertices are already colored, and $y_j$ is the only vertex of color $c_5$. However, it can be recolored with $c_3$, so it is not an acyclic b-vertex, a contradiction: $c$ contains no acyclic b-vertex that fulfills the conditions of Case 1.

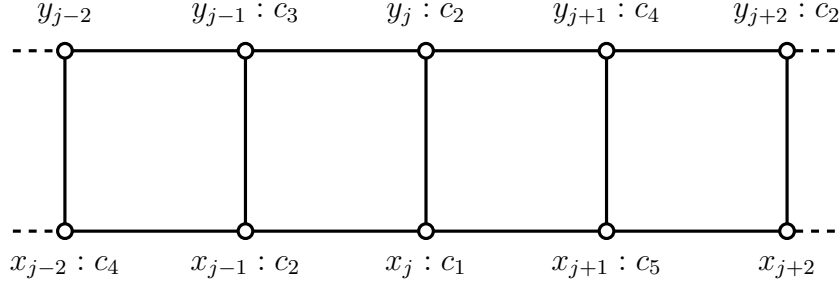
\begin{figure}[ht!]
\begin{center}

\begin{tikzpicture}
  [scale=0.8,auto=left,every node/.style={shape=circle,minimum size = 1pt, very thick}]
  \node[draw=black,scale=.5] (xj) at (0,0) {};
  \node[draw=black,scale=.5] (xj-1) at (-3,0) {};
  \node[draw=black,scale=.5] (xj-2) at (-6,0) {};
  \node[draw=white,scale=.5] (xj-3) at (-7,0) {};
  \node[draw=black,scale=.5] (xj+1) at (3,0) {};
  \node[draw=black,scale=.5] (xj+2) at (6,0) {};
  \node[draw=white,scale=.5] (xj+3) at (7,0) {};
  
  \node[draw=black,scale=.5] (yj) at (0,3) {};
  \node[draw=black,scale=.5] (yj-1) at (-3,3) {};
  \node[draw=black,scale=.5] (yj-2) at (-6,3) {};
  \node[draw=white,scale=.5] (yj-3) at (-7,3) {};
  \node[draw=black,scale=.5] (yj+1) at (3,3) {};
  \node[draw=black,scale=.5] (yj+2) at (6,3) {};
  \node[draw=white,scale=.5] (yj+3) at (7,3) {};

  \draw (0,-.6) node {$x_j:c_1$};
  \draw (-3,-.6) node {$x_{j-1}:c_2$};
  \draw (-6,-.6) node {$x_{j-2}:c_4$};
  \draw (3,-.6) node {$x_{j+1}:c_5$};
  \draw (6,-.6) node {$x_{j+2}$};

   \draw (0,3.6) node {$y_j:c_2$};
  \draw (-3,3.6) node {$y_{j-1}:c_3$};
  \draw (-6,3.6) node {$y_{j-2}$};
  \draw (3,3.6) node {$y_{j+1}:c_4$};
  \draw (6,3.6) node {$y_{j+2}:c_2$};

  \foreach \from/\to in {xj-2/xj-1,xj-1/xj,xj/xj+1,xj+1/xj+2,yj-2/yj-1,yj-1/yj,yj/yj+1,yj+1/yj+2,xj-2/yj-2,xj-1/yj-1,xj/yj,xj+1/yj+1,xj+2/yj+2}
  \draw[-,very thick] (\from) -- (\to);

   \foreach \from/\to in {xj-3/xj-2,yj-3/yj-2,xj+2/xj+3,yj+2/yj+3}
  \draw[dashed,very thick] (\from) -- (\to);

  
\end{tikzpicture}

\end{center}
\caption{Prisms: $c(x_{j-1})=c_2=c(y_{j})\neq c(x_{j+1})=c_5$.}
\label{prismCase1}
\end{figure}

\noindent{}\textbf{Case 2.} $x_j$ is type B vertex with $c(x_{j-1})=c(y_{j})\neq c(x_{j+1})$ (see Figure \ref{prismCase1}).\\
Let $c(x_{j-1})=c(y_{j})=c_2$ and $c(x_{j+1})=c_5$. There must exist two internally disjoint paths of even length with common endpoints $x_{j-1}$ and $y_j$, one colored with $c_2$ and $c_3$ and the other one with $c_2$ and $c_4$. The only possible configuration (up to the permutation of colors $\{c_3,c_4\}$) is presented in Figure \ref{prismCase1} and we have $c(y_{j-1})=c_3$, $c(x_{j-2})=c_4=c(y_{j+1})$ and $c(y_{j+2})=c_2$. Note that the internal vertices of the path colored with $c_2$ and $c_4$ are precisely $x_i$, $y_i$ or both ($x_i$ and $y_i$) for all $i\notin\{j-1,j,j+1\}$. 

Before we prove that this case cannot occur, let us describe the last option.

\begin{figure}[ht!]
\begin{center}

\begin{tikzpicture}
  [scale=0.8,auto=left,every node/.style={shape=circle,minimum size = 1pt, very thick}]
  \node[draw=black,scale=.5] (xj) at (0,0) {};
  \node[draw=black,scale=.5] (xj-1) at (-3,0) {};
  \node[draw=black,scale=.5] (xj-2) at (-6,0) {};
  \node[draw=white,scale=.5] (xj-3) at (-7,0) {};
  \node[draw=black,scale=.5] (xj+1) at (3,0) {};
  \node[draw=black,scale=.5] (xj+2) at (6,0) {};
  \node[draw=white,scale=.5] (xj+3) at (7,0) {};
  
  \node[draw=black,scale=.5] (yj) at (0,3) {};
  \node[draw=black,scale=.5] (yj-1) at (-3,3) {};
  \node[draw=black,scale=.5] (yj-2) at (-6,3) {};
  \node[draw=white,scale=.5] (yj-3) at (-7,3) {};
  \node[draw=black,scale=.5] (yj+1) at (3,3) {};
  \node[draw=black,scale=.5] (yj+2) at (6,3) {};
  \node[draw=white,scale=.5] (yj+3) at (7,3) {};
  
  \draw (0,-.6) node {$x_j:c_1$};
  \draw (-3,-.6) node {$x_{j-1}:c_2$};
  \draw (-6,-.6) node {$x_{j-2}:c_5$};
  \draw (3,-.6) node {$x_{j+1}:c_2$};
  \draw (6,-.6) node {$x_{j+2}:c_5$};

   \draw (0,3.6) node {$y_j:c_2$};
  \draw (-3,3.6) node {$y_{j-1}:c_4$};
  \draw (-6,3.6) node {$y_{j-2}$};
  \draw (3,3.6) node {$y_{j+1}:c_3$};
  \draw (6,3.6) node {$y_{j+2}$}; 
    
  \foreach \from/\to in {xj-2/xj-1,xj-1/xj,xj/xj+1,xj+1/xj+2,yj-2/yj-1,yj-1/yj,yj/yj+1,yj+1/yj+2,xj-2/yj-2,xj-1/yj-1,xj/yj,xj+1/yj+1,xj+2/yj+2}
  \draw[-,very thick] (\from) -- (\to);

   \foreach \from/\to in {xj-3/xj-2,yj-3/yj-2,xj+2/xj+3,yj+2/yj+3}
  \draw[dashed,very thick] (\from) -- (\to);

  
\end{tikzpicture}

\end{center}
\caption{Prisms: case $c(x_{j-1})=c(y_{j})=c(x_{j+1})=c_2$.}
\label{prismCase3}
\end{figure}
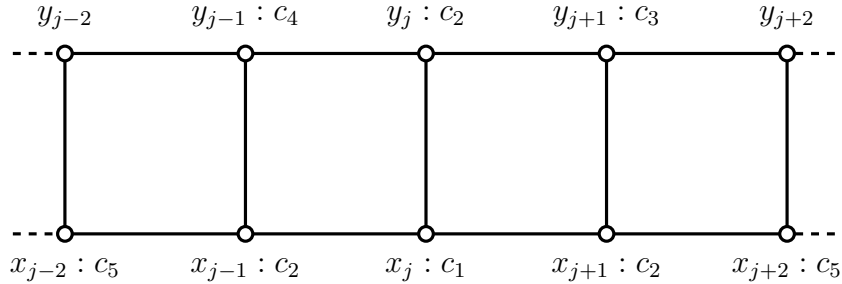

\noindent{}\textbf{Case 3.} $x_j$ is type A vertex with $c(x_{j-1})=c(x_{j+1})= c(y_{j})$ (see Figure \ref{prismCase3}).\\
Let $c(x_{j-1})=c(x_{j+1})=c(y_{j})=c_2$. There must exist three internally disjoint paths with endpoints in the set $\{x_{j-1},x_{j+1},y_j\}$, each of them colored with $c_2$ and one of the colors from $\{c_3,c_4,c_5\}$. The only possible configuration (up to the permutation of colors $\{c_3,c_4,c_5\}$) is presented in Figure \ref{prismCase3}. Note that in this case, the internal vertices of the path $P$ colored with $c_2$ and $c_5$ contain $x_i$ or $y_i$ for every $i\not\in\{j-1,j,j+1\}$.

Now we are ready to show that Cases 2 and 3 cannot occur. Note that in both cases there exists a bicolored path with both ends in the neighborhood of $x_j$ and the internal vertices being $x_i$ or $y_i$ for all $i\not\in\{j-1,j,j+1\}$. 

Assume that $x_j$ satisfies the conditions of Case 2. We will show that there exists no acyclic b-vertex of color $c_5$. Assume first that $x_{j+1}$ is such a vertex. In such a case, it needs to have at most two colors in its neighborhood, so $c(x_{j+2})=c_3$ is impossible. If $c(x_{j+2})=c_1$, then $x_{j+1}$ would satisfy the properties from Case 1, which is not possible, as we already know. If $c(x_{j+2})=c_4$, $x_{j+1}$ could satisfy the properties from Case 2. But this is impossible because $c(y_j)=c_2=c(y_{j+2})$. Hence, $x_{j+1}$ is not an acyclic b-vertex of color $c_5$.

Assume next that $x_i$ or $y_i$ is an acyclic b-vertex of color $c_5$ for some $i\not\in \{j-2, j-1, j, j+1, j+2 \}$. This would mean that there is a path colored with two colors, none of which is $c_5$, with both ends in the neighborhood of $x_i$ ($y_i$, respectively) having among its internal vertices in particular some of $x_{j-1},x_j,x_{j+1},y_{j-1},y_j,y_{j+1}$, which is impossible. The contradiction shows that Case 2 cannot occur.

Assume now that $x_j$ satisfies the conditions of Case 3. We will show that there exists no acyclic b-vertex of color $c_2$. Clearly, $x_{j-1}$, $x_{j+1}$, and $y_j$ are not such vertices because they have three different colors in their neighborhood. Thus, if $x_i$ or $y_i$ is an acyclic b-vertex of color $c_2$, then $i\not\in\{j-2,j-1,j,j+1,j+2\}$ ($i\not\in\{j-1,j,j+1\}$, respectively). However, just like in Case 2, this would imply the existence of a respective path colored with two colors, none of which is $c_2$, having among its internal vertices, in particular at least one of $x_{j-1}$, $x_{j+1}$, and $y_j$, a contradiction. This means that Case 3 is impossible.

Since the same configurations are the only ones possible for every acyclic b-vertex (up to the permutations of colors), we conclude that there exists no acyclic b-coloring with five colors, and we have $A_b(G(n,1))<5$. By Theorem \ref{bacycliccubic} it follows that $A_b(G(n,1))=4$. \qed


\section{$(0,j)$-prisms}

 As in \cite{ref_CFK2}, we generalize prisms.  Let the $(0,j)$-\emph{prism} of order $2n$, where $n$ and $j$ are even, be the graph with two vertex-disjoint cycles
$R_n^i=v_0^i,\dots,v_{n-1}^i$ for $i\in\{1,2\}$ of length $n$ called \emph{rims} and
edges $v_0^1v_0^2,v_2^1v_2^2,v_4^1v_4^2,\dots$ and
$v_1^1v_{j+1}^2,$ $v_3^1v_{3+j}^2,v_5^1v_{5+j}^2,\dots$
called \emph{spokes} of \emph{type} $0$ and \emph{type} $j,$ respectively.
It is easy to observe that $(0,j)$-prism is a cubic graph and is isomorphic to $(0,-j)$-prism, $(j,0)$-prism and $(-j,0)$-prism. We can therefore always assume that $0\leq j\leq\frac{n}{2}.$
In our terminology, the usual prism is referred to as a $(0,0)$-prism. We will denote $(0,j)$-prism with $n=2m$ vertices and $3m$ edges by \projn.

\begin{theorem}
If $j>0$ and $n\geq 5(j+2)$, then $A_b(Pr_n(0,j))=5$.
\end{theorem}

\begin{figure}[ht!]
\begin{center}
\begin{tikzpicture}[scale=1,style=very thick,x=1cm,y=1cm]
\def\vr{3pt} 
\path (0,0) coordinate (a);
\path (2,0) coordinate (b);
\path (4,0) coordinate (c);

\path (7,0) coordinate (e);
\path (9,0) coordinate (f);
\path (11,0) coordinate (g);

\path (-2,2) coordinate (a0);
\path (0,2) coordinate (a1);
\path (2,2) coordinate (b1);
\path (4,2) coordinate (c1);

\path (7,2) coordinate (e1);
\path (9,2) coordinate (f1);
\path (11,2) coordinate (g1);

\draw  (b) -- (c);
\draw (e) -- (f) -- (g) -- (g1) -- (f1) -- (e1);
\draw (c1) -- (b1) -- (a1) -- (a0);
\draw (b) -- (b1);
\draw (a1) -- (f);
 \draw[dashed,very thick] (c) -- (e);
 \draw[dashed,very thick] (c1) -- (e1);

\draw (b) [fill=white] circle (\vr);
\draw (c) [fill=white] circle (\vr);
\draw (e) [fill=white] circle (\vr);
\draw (f) [fill=white] circle (\vr);
\draw (g) [fill=white] circle (\vr);
\draw (a0) [fill=white] circle (\vr);
\draw (a1) [fill=black] circle (\vr);
\draw (b1) [fill=white] circle (\vr);
\draw (c1) [fill=white] circle (\vr);

\draw (e1) [fill=white] circle (\vr);
\draw (f1) [fill=white] circle (\vr);
\draw (g1) [fill=white] circle (\vr);

\draw[anchor = north] (b) node {$v_{2i+2}^{2,k}:c_{k}^1$};
\draw[anchor = north] (c) node {$v_{2i+3}^{2,k}:c_{k}^3$};
\draw[anchor = north] (e) node {$v_{2i+j}^{2,k}:c_{k}^1$};
\draw[anchor = north] (f) node {$v_{2i+1+j}^{2,k}:c_{k}^3$};
\draw (11.5,-0.43) node  {$v_{2i+2+j}^{2,k}:c_{k}^2$};
\draw[anchor = south] (a0) node {$v_{2i}^{1,k}:c_{k}^4$};
\draw[anchor = south] (a1) node {$v_{2i+1}^{1,k}:c_{k}^0$};
\draw[anchor = south] (b1) node {$v_{2i+2}^{1,k}:c_{k}^3$};
\draw[anchor = south] (c1) node {$v_{2i+3}^{1,k}:c_{k}^2$};
\draw[anchor = south] (e1) node {$v_{2i+j}^{1,k}:c_{k}^3$};
\draw[anchor = south] (f1) node {$v_{2i+1+j}^{1,k}:c_{k}^2$};
\draw (11.5,2.43) node  {$v_{2i+2+j}^{1,k}:c_{k}^3$};
\end{tikzpicture}
\end{center}
\caption{An acyclic b-vertex $v_{2i+1}^{1,k}$ for \projn.}
\label{projn}
\end{figure}
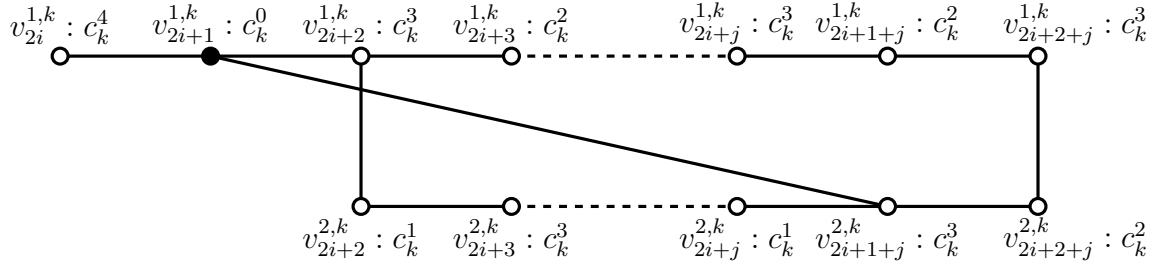

\noindent {\textbf{Proof.}} 
As in the proof of Theorem~\ref{thm_petersen}, we present the respective acyclic b-coloring with five colors. No vertex can be a b-vertex, and we construct an acyclic b-coloring with acyclic b-vertices of type B in every color. Consider the vertex $v_{2i+1}^1$ for some $0\leq i\leq n-1$ and the two cycles  $$C_k^1=v_{2i+1}^{1,k}v_{2i+2}^{1,k}v_{2i+2}^{2,k}v_{2i+3}^{2,k}\ldots v_{2i+1+j}^{2,k}v_{2i+1}^{1,k}\text{ and }C_{k}^2=v_{2i+1}^{1,k}v_{2i+2}^{1,k}\ldots v_{2i+2+j}^{1,k} v_{2i+2+j}^{2,k}v_{2i+1+j}^{2,k}v_{2i+1}^{1,k}.$$ 
The colors assigned to vertices of $C_k^1$ and $C_k^2$ are $c(v_{2i+1}^{1,k})=c_{k}^0$, $c(v_{2i+2}^{1,k})=c(v_{2i+4}^{1,k})=\cdots=c(v_{2i+2+j}^{1,k})=c(v_{2i+3}^{2,k})=c(v_{2i+5}^{2,k})=\cdots=c(v_{2i+1+j}^{2,k})=c_{k}^3$, $c(v_{2i+2}^{2,k})=c(v_{2i+4}^{2,k})=\cdots=c(v_{2i+j}^{2,k})=c_{k}^1$, $c(v_{2i+3}^{1,k})=c(v_{2i+5}^{1,k})=\cdots=c(v_{2i+1+j}^{1,k})=c(v_{2i+2+j}^{2,k})=c_{k}^2$ and $c(v_{2i}^{1,k})=c_{k}^4$, where $c_k^0$, $c_{k}^1$, $c_{k}^2$, $c_{k}^3$ and $c_{k}^4$ are pairwise distinct. Now, $v_{2i+1}^{1,k}$ is an acyclic b-vertex for color $c_{k}^0$, since colors $c_k^3$ and $c_k^4$ belongs to its neighbors, while $C_k^1$ is a $(c_k^1)_{v_{2i+1}^{1,k}}$-cycle and $C_k^2$ is $(c_k^2)_{v_{2i+1}^{1,k}}$-cycle. The situation is illustrated in Figure \ref{projn}.

One can repeat the presented configuration by choosing $i=\frac{k(j+2)}{2}$ for $0\leq k\leq 4$. This identifies the vertices $v_{2i+2+j}^{1,k}$  with $v_{2i}^{1,k+1}$, which can be done, since the set of mappings $c_k^3\rightarrow c_k^4$ defined for $0\leq k\leq 4$ can be obviously extended to the permutations giving the correct colorings of consecutive segments of graph even if one cyclically joins them. For the remaining vertices, we can apply a greedy procedure to obtain a proper coloring with five colors. Similarly, as in the proof of Theorem~\ref{acycliccubic}, we can avoid 2-colored cycles.~\qed




\section{Final remarks}

The ultimate goal for a cubic graph $G$ is to decide whether $A_b(G)=4$ or $A_b(G)=5$. We have presented some tools in this direction, but this is generally still open. 

\begin{problem}
What is the computational complexity of deciding whether $A_b(G)=4$ or $A_b(G)=5$ for any cubic graph $G$?
\end{problem}

Note that we do not know this for a general graph, see Problem 6.5 from \cite{AnCiPe}. While this may be a difficult task in the general case, we can find partial answers for certain classes of cubic graphs, for instance, snarks.

\begin{problem}
What is $A_b(G)$ of a snark?
\end{problem}

Some ideas on how to show that $A_b(G)=4$ for a cubic graph $G$ can be observed in Theorem \ref{notupper}. The reason there is $A_b(G)=4$ is that every vertex $v$ has acyclic degree $d_G^a(v)$ which equals three as well (for a definition of acyclic degree see Section 4 of \cite{AnCiPe}.) So, we can call such graphs acyclically cubic graphs and for them, we know that $A_b(G)=4$ (whenever $G$ is not $K_3\Box K_2$). 

\begin{problem}
Find classes/families of acyclically cubic graphs.
\end{problem}

We conclude with a conjecture that applies to arbitrary graphs, not only cubic ones. It seems that the process from Theorem~\ref{acycliccubic} can be extended for any graph $G$ such that a b-coloring is also acyclic if the girth $g(G)$ is big enough. 

\begin{conjecture}
Let $G$ be a graph such that $g(G)> 2\varphi(G)$. Then $A_b(G)\geq \varphi(G).$
\end{conjecture}


\end{document}